\documentclass[11pt]{article}
\usepackage{amssymb}
\usepackage{amsmath}
\usepackage{array}

\DeclareMathOperator{\RRe}{Re} \DeclareMathOperator{\IIm}{Im}

\multlinegap=0em \DeclareFontFamily{T1}{msb}{}
\DeclareFontShape{T1}{msb}{m}{ol}{<5> <6> <7> <8> <9> gen * msbm
<10> <10.95> <12> <14.4> <17.28> <20.74> <24.88> msbm10}{}
\DeclareSymbolFont{AMSb}{T1}{msb}{m}{ol} \multlinegap=0em

\renewcommand{\S}{\mathhexbox278}
\renewcommand{\le}{\operatorname{\leqslant}}
\renewcommand{\ge}{\operatorname{\geqslant}}

\def\hf{{\textstyle{\frac12}}}
\def\a{\alpha}\def\b{\beta}
\def\d{{\,\rm d}}
\def\e{\varepsilon}

\def\G{\Gamma} \def\g{\gamma}

\def\s{\sigma}
\def\z{\zeta}

\def\={\;=\;}

\def\zt{\zeta(\hf+it)}

\begin{document}

\begin{center} \large\bf  On the distribution of  values of the argument
of the Riemann zeta-function
\end{center}

\medskip
\begin{center}
\bf Aleksandar P. Ivi\'c and Maxim A. Korolev
\end{center}
\medskip

{\small {\bf Abstract:} Let $S(t) \;:=\; \frac{\displaystyle 1}{\displaystyle \pi}\arg \z(\hf + it)$. We prove that,
for $T^{\,27/82+\e} \le H \le T$, we have
$$
{\rm mes}\Bigl\{t\in [T, T+H]\;:\; S(t)>0\Bigr\} = \frac{H}{2} + O\left(\frac{H\log_3T}{\e\sqrt{\log_2T}}\right),
$$
where the $O$-constant is absolute.
A similar formula holds for the measure of the set with $S(t)<0$, where $\log_kT = \log(\log_{k-1}T)$.
This result is derived from an asymptotic formula
for the distribution of values of $S(t)$, which is uniform in the relevant parameters,
and this  is of crucial importance. This in fact depends on the distribution of values of the
Dirichlet polynomial which approximates $S(t)$, namely ($p$ denotes primes)
$$V_{y}(t)\,=\,\sum\limits_{p\le y}\frac{\sin{(t\log{p})}}{\sqrt{p}}.$$}


\vspace{0.3cm}


\medskip
\textbf{1. Introduction}

\medskip
The Riemann zeta-function
$$
\z(s) = \sum_{n=1}^\infty n^{-s}\qquad(\Re s >1),
$$
derives a part of its importance through the fact that it is a function
of the complex variable $s = \s+it\;(\s, t \in \mathbb R)$. As such, it possesses
meromorphic continuation to $\mathbb C$, its only singularity being the simple pole at $s=1$.
For an extensive account on $\z(s)$ the reader is referred to the monographs of
E.C.~Titchmarsh \cite{Titchmarsh} and the first author \cite{Ivic1}.

\medskip
Nevertheless, there are several important functions of the real variable $t$ connected to
the zeta-function. They contain much relevant information concerning $\z(s)$.
 One of these functions is Hardy's function
$$
Z(t)  := \zt\bigl(\chi(\hf+it)\bigr)^{-1/2},
$$
where $\z(s) \;=\; \chi(s)\z(1-s)$ is the functional equation for the zeta-function,
so that
$$
\chi(s) = \frac{\G(\hf(1-s))}{\G(\hf s)}\pi^{s-1/2}.
$$
It follows that $Z(t)$ is a smooth, real-valued function of the real variable
$t$, for which $$|Z(t)|=|\zt|.$$ Therefore the  real zeros of $Z(t)$ correspond to the
complex zeros of $\z(s)$ of the form $\hf + it$, namely to the zeros on the ``critical
line'' $\s = \hf$. For an extensive account  of $Z(t)$, see the first author's monograph \cite{Ivic2}.

\medskip
The second important function, which is the main subject of the present paper, is
the argument function
$$
S(t) \;:=\; \frac{1}{\pi}\arg \z(\hf + it)\qquad(t>0, \;t \ne \g),
$$
where $\rho = \b+i\g$ denotes generic complex zeros of $\z(s)$. If $t=\g$, we follow
Titchmarsh \cite{Titchmarsh} and define \footnote{A.~Selberg \cite{Selberg_1944} defines
$$
S(t) = \hf(S(t+0) + S(t-0))\qquad(t=\g).
$$
For our purposes it does not matter which definition of $S(t)$ is used.}
$$
S(t) \;=\; S(t+0)\qquad(t=\g).
$$
Unlike $Z(t)$, the argument function $S(t)$ is not continuous, but has jumps (discontinuities
of the first kind) at ordinates $\g$ of zeta-zeros. On every interval $(\g, \g^+)$, where
$\g, \g^+\;(0< \g < \g^+)$ are consecutive ordinates of zeta-zeros, the
function $S(t)$ is monotonically decreasing and
$$
S'(t) = - \frac{1}{2\pi}\log \frac{t}{2\pi} + O\left(\frac{1}{t^2}\right),
\;
S''(t) = - \frac{1}{2\pi t} + O\left(\frac{1}{t^3}\right).
$$
This follows from the Riemann-von Mangoldt formula
$$
N(t) = \sum_{0<\g\le t}1 = \frac{t}{2\pi}\log \frac{t}{2\pi} - \frac{t}{2\pi} + \frac{7}{8}
+ S(t) + f(t),
$$
where multiple zeros are counted with their multiplicities and $f(t)$ is a smooth function such that
$$
f(t) = O\left(\frac{1}{t}\right), \quad f'(t) = O\left(\frac{1}{t^2}\right),
\quad f''(t) = O\left(\frac{1}{t^3}\right).
$$
One has the bounds
$$
S(T) = O(\log T),\qquad S(T) = O\left(\frac{\log T}{\log\log T}\right)\quad(\rm{RH}),\eqno(1.1)
$$
where the first one is unconditional, and (RH) means that the bound in question
holds under the Riemann Hypothesis (all complex zeros of
$\z(s)$ have real parts equal to 1/2). Note that the RH gives an improvement of only the factor of
$\log\log T$ over the unconditional bound $S(T) = O(\log T)$.

\bigskip
With real-valued functions such as $Z(t)$ and $S(t)$, which take on positive and negative values that
are not regularly distributed, one may naturally ask: what is the measure of the sets in $[T, T+H]$
where these functions are positive or negative? Here $0<H = H(T)\le T$ is to be suitably fixed.

\medskip
This problem for $Z(t)$ was recently investigated by S.M.~Gonek and the first author \cite{GonekIvic}. It was
proved there ($\rm{mes}\{\cdot\}$ denotes measure) that
$$
{\rm{mes}}\bigl\{T\in [T,2T]\;:\; Z(t)>0\bigr\} \;\gg\; T, \quad {\rm{mes}}\bigl\{T\in [T,2T]\;:\; Z(t)<0\bigr\}\; \gg\; T.
 \eqno(1.2)
$$
Of course, one of the bounds in (1.2) must hold, since ${\rm{mes}}\bigl\{\,[T, 2T]\,\bigr\} = T$, but {\it a priori} one
cannot say which one holds. If one assumes the RH and H.L.~Montgomery's pair correlation conjecture
(see \cite{Montgomery}), then (1.2) can be improved to
$$
{\rm{mes}}\bigl\{T\in [T,2T]\;:\; Z(t)>0\bigr\} \ge 0.32909T, \quad {\rm{mes}}\bigl\{T\in [T,2T]\;:\; Z(t)<0\bigr\} \ge 0.32909T,
$$
provided that $T$ is sufficiently large. The pair correlation conjecture states that, if one assumes the RH,
$$
\sum_{0<\g,\g'\le T,\frac{2\pi\a}{\log T}\le \g'-\g\le\frac{2\pi\b}{\log T}}1\;
\sim \;\left(\int\limits_\a^\b\left[1-\left(\frac{\sin(\pi u)}{\pi u}\right)^2\right]\d u
+ \delta(\a,\b)\right)\frac{T}{2\pi}\log T
$$
as $T\to\infty$. Here $\g, \g'$ denote arbitrary ordinates of zeta-zeros,  $\a < \b$ are fixed numbers,
and $\delta(\a,\b) =1$ if $0\in[\a,\b]$ and $\delta(\a,\b) =0$ otherwise.

\medskip
As to the distribution of positive and negative values of $S(t)$, there appear to be no results in
the literature so far. However, since this is not a continuous function, it makes also sense to
investigate $M(T)$, the number of sign changes of $S(t)$ in $(0, T]$. A.~Selberg was the first to
obtain significant results concerning this problem, and he proved (see \cite{Selberg_1944})  that
$$
M(T+H) - M(T) \ge H(\log T)^{1/3}\exp\bigl(-C\sqrt{\log\log T}\,\bigr)\quad(T\ge T_0(\e), C>0, H = T^{1/2+\e}).
$$
Selberg's methods were further developed in the thesis of K.-M.~Tsang \cite{Tsang_1984}.
For basic results on $M(T)$, and on $S(T)$ in general, the reader is referred to the surveys \cite{KarKor_2005} and \cite{KarKor_2006} of A.A.~Karatsuba and the second author.

\medskip
\textbf{2. Statement of results}

\medskip
Our main aim is to investigate the distribution of positive and negative values of $S(t)$.

\medskip
THEOREM 1. {\it Suppose that $0<\varepsilon<10^{-3}$ is an arbitrary small fixed constant,
$T\ge T_{0}(\varepsilon)$, $T^{c+\e} \le H \le T$,
where $c = \tfrac{27}{82}$. Then,  for any real $a$ and $b$, $a<b$,}
\begin{equation*}
{\rm mes}\Biggl\{t\in[T, T+H] \,:\, a < \frac{\pi\sqrt{2}\,S(t)}{\sqrt{\log\log T\mathstrut}} \le b\Biggr\}
= \frac{H}{\sqrt{2\pi}}\left(\int_a^b e^{-v^2/2}\d v + O\left(\frac{\log_3T}{\varepsilon\sqrt{\log_2T}}\right)\right),\eqno(2.1)
\end{equation*}
{\it where the O-constant is absolute. }

\medskip
The reason that on the left-hand side of (2.1) one has $\pi\sqrt{2}\,S(t)/{\sqrt{\log\log T\mathstrut}}$ is that this function equals unity in the mean square sense.
This was established first by A.~Selberg \cite{Selberg_1946}, and for short intervals
we refer to a result of A.A.~Karatsuba \cite{Karatsuba_1996}.
Namely, he proved (see his Theorem B with $k=1$) that
$$
\int_T^{T+H}S^2(t)\d t = \frac{H}{(\pi\sqrt{2})^{2\mathstrut}}\log\log T + O\bigl(H\sqrt{\log\log T\mathstrut}\,\bigr)
$$
for $H = T^{\,27/82+\e}$.

\medskip
The assertion (2.1) was obtained first by K.-M.~Tsang (see \cite{Tsang_1984}, Theorem 6.1),
where he used Selberg's density theorem \cite{Selberg_1946} for the zeros of the
Riemann zeta-function $\zeta(s)$ lying in the rectangle $\sigma<\RRe s\le 1$,
 $T<\IIm s\le T+H$,  where $T^{\,1/2+\varepsilon}\le H\le T$. A.A.~Karatsuba \cite{Karatsuba_1996}
 obtained a density theorem which enabled
R.N.~Boyarinov \cite{Boyarinov_2011} to transfer the analogues of the main theorems from
 \cite{Tsang_1984} to the case $T^{\,c+\varepsilon}\le H\le T$, $c = \tfrac{27}{82}$
 (his paper \cite{Boyarinov_2011} contains only formulations of the assertions).
 For our purposes, the uniformity in the parameters $a$ and $b$ is crucial, since in our applications $b$ will be taken to be of the order of $(\log T)/\sqrt{\log\log{T\mathstrut }}$.
In Section 4 we shall give a full proof of  formula (2.1) in order to underline the fact the
 bound for the remainder term in (2.1) does not depend on the parameters $a,b\in \mathbb{R}$.

\medskip
THEOREM 2.  {\it Suppose that $0<\varepsilon<10^{-3}$ is an arbitrary small fixed constant,
$T\ge T_{0}(\varepsilon)$, $T^{c+\e} \le H \le T$,
where $c = \tfrac{27}{82}$. Then}
$$
{\rm{mes}}\Bigl\{t\in [T, T+H]\;:\; S(t)>0\Bigr\} = \frac{H}{2} + O\left(\frac{H\log_3T}{\e\sqrt{\log_2T}}\right),\eqno(2.2)
$$
{\it where the O-constant is absolute}.

\medskip
Since ${\rm{mes}}\bigl\{[T, T+H]\bigr\} = H$, it follows immediately that (2.2) implies also
$$
{\rm{mes}}\Bigl\{t\in [T, T+H]\;:\; S(t)<0\Bigr\} = \frac{H}{2} + O\left(\frac{H\log_3T}{\e\sqrt{\log_2T}}\right).\eqno(2.3)
$$
The lower bound $T^{\,27/82+\e}$ for $H$ is certainly not optimal, but just a convenient one that the method allows.
This bound originated in the work of Karatsuba \cite{Karatsuba_1984} on the zeros of $\z(s)$ on short intervals of the critical line.
In his subsequent works Karatsuba (see e.g., \cite{Karatsuba_1996}) used the same method.
Indeed, as remarked in the review Zbl.0545.10026, if instead of (35) on p. 580 of \cite{Karatsuba_1984} one uses
$$
\mathbb{C}(u) \;\ll\; T^\kappa|h+\xi_2-\xi_1|M^{\lambda-2\kappa},
$$
where $(\kappa, \lambda)$ is an exponent pair (see e.g., Chapter 2 of \cite{Ivic1}), one obtains that the constant $\tfrac{27}{82}+\e$
in Karatsuba's works can be replaced by $\a+\e$, where
$$
\a := \inf_{(\kappa,\lambda)}(\kappa + \lambda - \hf) = 0.329021\ldots\, (<  \tfrac{27}{82} = 0.329286\ldots\,),
$$
the so-called Rankin's constant in the theory of exponent pairs. Further small improvements can be attained if,
instead of the classical exponent pairs, one uses the new exponent pairs obtained by M.N.~Huxley and recently
by J.~Bourgain \cite{Bourgain}. In the context of the present work, the constant $\tfrac{27}{82}$ appears via Lemma 12, which
ultimately depends on the arguments of Karatsuba's paper (op. cit.), hence the above improvement of the constant $\tfrac{27}{82}$
holds also in our case. In particular, by taking $c$ slightly less than $\tfrac{27}{82}$, we shall
obtain Theorem 1 and Theorem 2 with $c = \tfrac{27}{82}$, but without $\e$.

It may be true that even $T^\e$ could be taken in (2.1)--(2.3), but this
conjecture is certainly out of reach of the present methods.
Note that (2.2) and (2.3) are the analogue of (1.2) for $S(t)$, only they represent  true asymptotic formulas, and moreover
they hold over the ``short'' interval $[T, T +H]$. The asymptotic formula (2.2) is an easy consequence of (2.1), but
it is stated as a theorem because it is attractive and represents an analogue of (1.2) for the function $S(t)$.
To see how (2.2) follows from (2.1), set $a = 0, b = \pi C\sqrt{2}(\log T)/\sqrt{\log\log T\mathstrut}$, where $C > 0$ is such a constant that $|S(T)| \le C\log T$
for $T\ge T_0$. Such a constant must exist in view of the upper bound (1.1). Then (2.2) follows, since
$$
\int_0^\infty e^{-v^2/2}\d v = \sqrt{\frac{\pi}{2}}.
$$

\medskip

\textsc{THEOREM 3.} \emph{Suppose that $1 < H_0 < H $ and}
\begin{equation*}
y_{0}<y\le\exp{\biggl(\frac{\log{H}}{c\log\log{H}}\biggr)},\quad c = 1.4\cdot 10^{4}.\eqno(2.4)
\end{equation*}
\emph{Then, for any real $T$, $\alpha$ and for any integer $\nu$ such that}
\begin{equation*}
0\,\le\,\nu\,\le\,\exp{\bigl((\log{y})^{1/4}\bigr)}\eqno(2.5)
\end{equation*}
\emph{the following relation holds}:
\begin{equation*}\begin{split}&
\int_{T}^{T+H}\text{sgn}\bigl(V_{y}(t)-\alpha\bigr)\d t \\&
=\,\frac{H}{\sqrt{2\pi}}\sum\limits_{n = 0}^{\nu}\frac{\Phi_{n}}{(2\sigma)^{n}}
\int_{-\infty}^{+\infty}e^{-\,u^{2}/2}H_{2n}(u)\,\text{sgn}\biggl(u\sqrt{\frac{\sigma}{2}}\,-\,\alpha\biggr)\d u\,+\,O\bigl(H\Delta\bigr),
\quad\quad\quad\qquad\qquad\quad(2.6)\end{split}
\end{equation*}
\emph{where $V_y(t)$ and $\s$ are defined by} (3.1), $\Phi_n$ \emph {before Lemma 4, Hermite polynomial $H_{k}(x)$ by} (3.3),
\[
\Delta\,=\,\frac{{(\log{y})}^{-1/2}}{\log\log{y}}\,+\,\biggl(\frac{\log{y}}{\log{H}}\biggr)^{\! 1/2}\,+\,
\frac{1}{\nu+1}\biggl(\frac{\log\log{(\nu+2)}+1}{\sigma}\biggr)^{\!\nu+1},\eqno(2.7)
\]
\emph{and the implied O-constant is absolute.}

\medskip
It is Theorem 3 which, in spite of its unwieldy formulation, is the deepest and most difficult of our results.
The crux of the matter is the statement at the and of the theorem that {\it the implied O-constant is absolute}.
This will allow us to pass from $V_y(t)$ to $S(t)$ and deduce Theorem 1. Formulas (2.6) and (2.7) do not seem
to have appeared before in any form.

\medskip
\textsc{Corollary 1. } \emph{Let $\alpha<\beta$. Then, under the conditions of Theorem 3, we have}
\[
\int_{T}^{T+H}\chi_{\alpha,\beta}\bigl(V_{y}(t)\bigr)\d t\,=\,\frac{H}{\sqrt{2\pi}}
\sum\limits_{n=0}^{\nu}\int_{-\infty}^{+\infty}e^{-u^{2}/2}H_{2n}(u)\chi_{a,b}(u)\d u\,+\,O\bigl(H\Delta\bigr),\eqno(2.8)
\]
\emph{where $a = \alpha\sqrt{2/\sigma}$, $b = \beta\sqrt{2/\sigma}$, $\chi_{p,q}(x)$ is the characteristic function of the segment
$[p,q]$ and the implied constant is absolute.}

\vspace{0.5cm}

\medskip
\textbf{3. Notation}

\vspace{0.5cm}
In this section we present, for the convenience of the reader, some standard notation which
will be used later in the body of the text.

\medskip
For real $t$ and for sufficiently large $y>y_{0}>2$ we set
\[
V_{y}(t)\,=\,\sum\limits_{p\le y}\frac{\sin{(t\log{p})}}{\sqrt{p}},\quad
\sigma\,= \sigma(y) \;=\; \,\sum\limits_{p\le y}\frac{1}{p},\eqno(3.1)
\]
where $p$ runs through prime numbers. Further, let $J_{0}(z)$ be the Bessel function of the first kind, that is,
\[
J_{0}(z)\,=\,\sum\limits_{n=0}^{+\infty}\frac{(-1)^{n}}{(n!)^{2}}\biggl(\frac{z}{2}\biggr)^{\! 2n}\qquad(\forall z \in \mathbb C).
\]
Next, for $\alpha<\beta$, we denote by $\chi_{\alpha,\beta}(x)$  the characteristic function of the segment $[\alpha,\beta]$, namely
\begin{equation*}
\chi_{\alpha,\beta}(x)\,=\,
\begin{cases}
1, & \text{for}\quad \alpha\le x\le \beta,\\
0, & \text{otherwise}.
\end{cases}
\end{equation*}
Also the sign function sgn$(x)$ is commonly defined as
\begin{equation*}
{\rm sgn}(x)\,=\,
\begin{cases}
1, & \text{for}\quad  x > 0,\\
0, & \text{for}\quad x = 0,\\
-1 &  \text{for}\quad x<0.
\end{cases}
\end{equation*}
As usual, $\pi(y)$ is the number of primes $p$ not exceeding $ y$, $\Omega(n)$ is the number
of prime divisors of $n$ counted with multiplicities, $\Omega(1) =0$,
$P(n)$ is the largest prime divisor of $n>1, P(1)=1$. By $\theta,\theta_{1},\ldots$ we denote complex numbers with modulus at most $1$.
By $B$ we denote the so-called Mertens constant, that is,
\[
B\,=\,\lim_{x\to +\infty}\biggl(\;\sum\limits_{p\le x}\frac{1}{p}\,-\,\log\log{x}\biggr)\,=\,0.26149\,72128\,\ldots\;.
\]
Therefore, by the prime number theorem,
\[
\sigma\,= \sigma(y) \;=\; \log\log y + B + O\bigl(\exp(-\sqrt{\log y}\,)\bigr).
\]

\medskip
Next, let
\[
K(x)\,=\,\biggl(\frac{\sin{(\pi x)}}{\pi x}\biggr)^{\! 2}.\eqno(3.2)
\]
By $\text{mes}\{{\cal E}\}$ and by $|{\cal E}|$ we denote the Lebesgue measure of the set ${\cal E}\subset \mathbb{R}$.
All the constants in the symbols $O$ and $\ll$ are absolute.

\medskip
For an integer $n\ge 0$, the Hermite polynomial $H_{n}(x)$ is defined by the relation
\[
H_{n}(x)\,=\,(-1)^{n}e^{\,x^{2}/2}\,\frac{\d^{n}}{\d x^{n}}\bigl(e^{-x^{2}/2}\bigr).
\]

\vspace{0.3cm}

One can show that $H_{0}(x)\equiv 1$, $H_{1}(x) = x$ and $H_{n}(x) = xH_{n-1}(x)-(n-1)H_{n-2}(x)$ for $n\ge 2$.
Hence,
\[
H_{2}(x)\,=\,x^{2}-1,\quad H_{3}(x)\,=\,x^{3}-3x,\quad H_{4}(x)\,=\,x^{4}-6x^{2}+3,\quad H_{5}(x)\,=\,x^{5}-10x^{2}+15x,
\]
and, more generally,
\begin{equation*}
H_{n}(x)\,=\,\sum\limits_{0\le j\le n/2}\frac{(-1)^{j}n!}{j!(n-2j)!}\,(2x)^{n-2j}.\eqno(3.3)
\end{equation*}

\vspace{0.5cm}

\textbf{4. The necessary lemmas}

\medskip
This section contains lemmas necessary for the proof of our theorems. Some of them are straightforward, but some
are elaborate and seem to be of independent interest.
\vspace{0.5cm}

\textsc{Lemma 1.} \emph{Let $n\ge 1$ be any integer and let}
\[
I_{n}\,=\,I_{n}(T;H)\,=\,\int_{T}^{T+H}V_y^{\,n}(t)\d t,
\]
\emph{where $V_y(t)$ is defined by} (3.1).\emph{ Then, for any $T$ and $H>1$, one has}
\begin{equation*}\label{lab-2-1}
I_{n}\,=\,H\delta_{n}\,+\,O\bigl((\tfrac{1}{2}\,ny)^{n/2}\bigr),\eqno(4.1)
\end{equation*}
\emph{where}
\begin{equation*}
\delta_{n}\,=\,
\begin{cases}
\frac{\displaystyle (2k)!}{\displaystyle k!\mathstrut}\,\frac{\displaystyle G^{(k)}(0)}{\displaystyle 2^{2k\mathstrut }}, & \textit{for even}\quad n = 2k,\\
0, & \textit{for odd}\quad n = 2k-1,
\end{cases}
\end{equation*}
\emph{and}
\[
G(z)\,=\,\prod\limits_{p\le y}J_{0}\biggl(2i\sqrt{\frac{z}{p}}\,\biggr).\eqno(4.2)
\]

\textsc{Proof.} All these relations were established in lemmas 3.3 and 3.4 of the thesis of K.-M.~Tsang \cite{Tsang_1984}.
The only difference is that the factor $\delta_{2k}$ in \cite{Tsang_1984} has the form
\[
\binom{2k}{k}\,\frac{S_{k}}{2^{2k\mathstrut}},\quad S_{k}\,=\,\sum\frac{1}{p_{1}\ldots p_{k}},
\]
where  summation is taken over all ordered $k$-tuples $(p_{1},\ldots,p_{k})$ and $(q_{1},\ldots,q_{k})$ of primes $p_{1},\ldots , q_{k}\le y$ such that
\begin{equation*}
p_{1}\ldots p_{k}\,=\,q_{1}\ldots q_{k}.\eqno(4.3)
\end{equation*}
Now we shall transform the sum $S_{k}$ by a procedure which is due  to M.~Radziwi{\l}{\l} \cite{Radziwill_2011}.

\medskip
Given such a $k$-tuple $(p_{1},\ldots,p_{k})$, denote by $r_{1},\ldots, r_{s}$ all its different components and
by $\alpha_{1},\ldots,\alpha_{s}$ their multiplicities. Hence, the equation (4.3) has
\[
\frac{(\alpha_{1}+\ldots+\alpha_{s})!}{\alpha_{1}!\ldots\alpha_{s}!}\,=\,\frac{k!}{\alpha_{1}!\ldots\alpha_{s}!}
\]
solutions in primes $q_{1},\ldots, q_{k}$.

Suppose now that all prime divisors of $n$  do not exceed $y$ and, moreover, $\Omega(n) = k$. Then the equation
$p_{1}\ldots p_{k} = n = q_{1}\ldots q_{k}$ has $(k!)^{2}(\alpha_{1}!\ldots\alpha_{s}!)^{-2}$ solutions
in primes $p_{1},\ldots, q_{k}\le y$. Writing $n$ in the canonical form $n = r_{1}^{\alpha_{1}}\ldots r_{s}^{\alpha_{s}}$
and summing over all such integers $n$, we represent $S_{k}$ as follows:
\[
S_{k}\,=\,\sum\limits_{\substack{\Omega(n) = k \\ P(n)\le y}}\frac{1}{n}\biggl(\frac{k!}{\alpha_{1}!\ldots\alpha_{s}!}\biggr)^{\!2}
\,=\,(k!)^{2}\sum\limits_{s\ge 1}\;\;\sum\limits_{\substack{\alpha_{1}+\ldots+\alpha_{s}=k \\ r_{1},\ldots,r_{s}\le y}}
\frac{(\alpha_{1}!)^{-2}}{r_{1}^{\alpha_{1}}}\ldots \frac{(\alpha_{s}!)^{-2}}{r_{s}^{\alpha_{s}}}.
\]
Define the multiplicative function $f(n)$ on prime powers $p^\a$ as follows:
\begin{equation*}
f(p^{\alpha})\,=\,
\begin{cases}
(\alpha!)^{-2}, & \text{if}\quad p\le y,\\
0, & \text{otherwise}.
\end{cases}
\end{equation*}
Then $S_{k} = (k!)^{2}\mathfrak{g}_{k}$, where
\[
\mathfrak{g}_{k}\,:=\,\sum\limits_{\substack{n=1 \\ \Omega(n) = k}}^{+\infty}\frac{f(n)}{n}.
\]
Now let us consider the (formal) power series
\[
G(z)\,=\,\sum\limits_{k = 0}^{+\infty}\mathfrak{g}_{k}\,z^{k}.
\]
Changing the order of summation and using the multiplicativity of the function $z^{\Omega(n)}$ over $n$, we find that
\begin{alignat*}{2}
G(z)&\,=\,\sum\limits_{k=0}^{+\infty}z^{k}\biggl(\;\sum\limits_{\substack{n=1 \\
\Omega(n) = k}}^{+\infty}\frac{f(n)}{n}\biggr)\,=\,\sum\limits_{n=1}^{+\infty}\frac{f(n)}{n}\,z^{\Omega(n)}\\&
=\,\prod\limits_{p}\biggl(1\,+\,\sum\limits_{\alpha = 1}^{+\infty}\frac{f(p^{\alpha})}{p^{\alpha}}\,z^{\alpha}\biggr)\,=\,
\prod\limits_{p\le y}\biggl(1\,+\,
\sum\limits_{\alpha = 1}^{+\infty}\frac{1}{(\alpha!)^{2}}\biggl(\frac{z}{p}
\biggr)^{\!\alpha}\biggr).
\tag{4.4}
\end{alignat*}
Each factor of the above product is a convergent series that converges on every compact domain of the complex plane.
Therefore, $G(z)$ is absolutely convergent series since it is a product of finite number of absolutely convergent series.
This justifies the change of  the order of summation in (4.4). Hence $G(z)$ is an entire function and
\[
\mathfrak{g}_{k}\,=\,\frac{G^{(k)}(0)}{k!},\quad S_{k}\,=\,k!\,G^{(k)}(0).
\]
Noting that
\[
1\,+\,
\sum\limits_{\alpha = 1}^{+\infty}\frac{1}{(\alpha!)^{2}}\biggl(\frac{z}{p}\biggr)^{\!\alpha}\,=\,
\sum\limits_{\alpha = 0}^{+\infty}\frac{(-1)^{\alpha}}{(\alpha!)^{2}}\biggl(\frac{2i}{2}\sqrt{\frac{z}{p}}\,
\biggr)^{\! 2\alpha}
\,=\,J_{0}\biggl(2i\sqrt{\frac{z}{p}}\,\biggr),
\]
we obtain the assertion of the lemma. $\Box$

\vspace{0.5cm}

Let  $\varphi(z) = e^{-z}J_{0}(2i\sqrt{z})$ and denote by $\varpi_{n}$, $n = 0,1,2,\ldots$ the coefficients of its expansion
into Taylor series:
\[
\varpi_{0}\,=\,1,\quad \varpi_{1}\,=\,0,\quad \varpi_{2}\,=\,-\frac{1}{4},\quad \varpi_{3}\,=\,\frac{1}{9},\quad
\varpi_{4}\,=\,-\,\frac{5}{192},\quad \ldots\,.
\]

\textsc{Lemma 2.} \emph{For any $n\ge 0$, one has} $|\varpi_{n}|< (n!)^{-1}$.

\vspace{0.5cm}

This is lemma 1.10 from \cite{Korolev_2015}.

\vspace{0.5cm}

\textsc{Lemma 3.} \emph{Suppose that $u$ is real and $k\ge 0$ is an integer. Then the following estimates hold}:

(a) $|G(-u^{2})|\le 1$ \emph{for any} $u$;

(b) $|G(-u^{2})|\le e^{-0.5\sigma u^{2}}$ \emph{for} $1\le u \le \tfrac{1}{4}e^{0.5(\log{y})^{1/3}}$;

(c) $0<G^{(k)}(0)\le \sigma^{k}$ \emph{for any} $k$.

\vspace{0.5cm}

\textsc{Proof.} The inequality (a) follows from the definition of $G(z)$ (see (4.2)) and the identity
\[
J_{0}(x)\,=\,\frac{1}{\pi}\int_{0}^{\pi}\cos{(x\sin{t})}\d t.
\]
Next, if $u$ satisfies the conditions of (b) then
\begin{equation*}\label{lab-2-4}
(4u)^{2}<e^{(\log{y})^{1/3}}<y.\eqno(4.5)
\end{equation*}
Using (4.5) together with (a) we get
\[
|G(-u^{2})|\,=\,\prod\limits_{p\le y}\biggl|J_{0}\biggl(-\,\frac{2u}{\sqrt{p}}\biggr)\biggr|\,\le\,
\prod\limits_{(4u)^{2}<p\le y}\biggl|J_{0}\biggl(-\,\frac{2u}{\sqrt{p}}\biggr)\biggr|.
\]
Given a prime $p$, $(4u)^{2}<p\le y$, we set $v = u/\sqrt{p}$. Hence, $|v|\le 1/4$, so, by Lemma 2 we have
\begin{multline*}
\biggl|J_{0}\biggl(-\,\frac{2u}{\sqrt{p}}\biggr)\biggr|\,=\,e^{-v^{2}}\bigl|e^{v^{2}}J_{0}\bigl(-2v\bigr)\bigr|\,=
\,e^{-v^{2}}\bigl|\varphi(v^{2})\bigr|\,\le\\
\le\,e^{-v^{2}}\sum\limits_{n=0}^{+\infty}|\varpi_{n}|v^{2n}\,\le\,e^{-v^{2}}\biggl(1\,+\,v^{2}\sum\limits_{n=2}^{+\infty}
\frac{1}{n! 16^{n-1}}\biggr)\,<\\
<\,e^{-v^{2}}\biggl(1\,+\,\frac{v^{2}}{31}\biggr)\,<\,e^{-v^{2}+v^{2}/31}\,=\,e^{-(30/31)v^{2}}.
\end{multline*}
Thus,
\[
|G(-u^{2})|\,<\,\prod\limits_{(4u)^{2}<p\le y}\exp{\biggl(-\,\frac{30}{31}\,\frac{u^{2}}{p}\biggr)}\,=\,
\exp{\biggl(-\frac{30}{31}\,u^{2}\,\sum\biggr)},\quad \sum\,=\,\sum\limits_{(4u)^{2}<p\le y}\frac{1}{p}.
\]
Next, we express $\Sigma$ as the difference $\sigma-\sigma'$, where
\[
\sigma = \sum\limits_{p\le y}\frac{1}{p}, \qquad\sigma'\,=\,\sum\limits_{p\le (4u)^{2}}\frac{1}{p}.
\]
By the inequality (3.20) from \cite{Rosser_Schoenfeld_1962}, we have
\begin{multline*}
\sigma'\,\le\,\log\log{(4u)^{2}}\,+\,B\,+\,\frac{1}{2\log^{2}{(4u)^{2}}}\,=\\
=\,\log\log{(4u)^{2}}\biggl(1\,+\,\frac{B}{\log\log{(4u)^{2}}\mathstrut}
+\frac{1}{(\log{(4u)^{2}})^{2}\log\log{(4u)^{2}}}\biggr)\,\le\\
\le\,\log\log{(4u)^{2}}\biggl(1\,+\,\frac{B}{\log\log{16}}\,+\,\frac{1}{2(\log{16})^{2}\log\log{16}}\biggr)
\,<\,\frac{4}{3}\,\log\log{(4u)^{2}}.
\end{multline*}
Using (4.5), we find that
\[
\log{(4u)^{2}}\,<\,(\log{y})^{1/3},\quad \log\log{(4u)^{2}}\,<\,\frac{1}{3}\log\log{y},
\]
and hence, by the inequality (3.19) from \cite{Rosser_Schoenfeld_1962}, we obtain
\[
\sigma'\,<\,\frac{4}{9}\,\log\log{y}\,<\,\frac{4}{9}\,\sigma,\quad \Sigma\,>\,\frac{5}{9}\,\sigma,\quad |G(-u^{2})|\,<\,\exp{\biggl(-\,\frac{30}{31}\cdot\frac{5}{9}\,\sigma u^{2}\biggr)}\,<\,e^{\,-0.5\sigma u^{2}}.
\]
Finally, to prove (c), we use the arguments from Lemma 1:
\[
G^{(k)}(0)\,=\,\frac{1}{k!}\sum\limits_{\substack{p_{1}\ldots p_{k} = q_{1}\ldots q_{k}\\ p_{1},\ldots, q_{k}\le y}}\frac{1}{p_{1}\ldots p_{k}}\,\le\,\frac{1}{k!}
\sum\limits_{p_{1},\ldots,p_{k}\le y}\frac{k!}{p_{1}\ldots p_{k}}\,=\,\sigma^{k}.
\]
The lemma is proved. $\Box$

\medskip

Consider the entire function
\[
\Phi(z)\,:=\,\prod\limits_{p\le y}\varphi\biggl(\frac{z}{p}\biggr).
\]
Then
\[
\Phi(z)\,=\,\prod\limits_{p\le y}J_{0}\biggl(2i\sqrt{\frac{z}{p}}\,\biggr)e^{\,-z/p}\,=\,e^{\,-\sigma z}G(z),
\]
so that the function $G(z)$ is expressed in the form
\begin{equation*}
G(z)\,=\,e^{\,\sigma z}\,\Phi(z).\eqno(4.6)
\end{equation*}
Denote by $\Phi_{n}$ the coefficients of the expansion of $\Phi(z)$ into Taylor series. In \cite[\S 1.2]{Korolev_2015}, one can find
the explicit expressions for the values $\Phi_{n}$, $0\le n\le 10$. In particular,
\begin{align*}
& \Phi_{0}\,=\,1,\quad \Phi_{1}\,=\,0,\quad \Phi_{2}\,=\,-\,\frac{1}{2^{2\mathstrut}}\,\sigma_{2},\quad
\Phi_{3}\,=\,\frac{1}{3^{2\mathstrut}}\,\sigma_{3},\\
& \Phi_{4}\,=\,-\,\frac{11}{2^{6\mathstrut}\!\cdot\!3}\,\sigma_{4}\,+\,\frac{1}{2^{5\mathstrut}}\,\sigma_{2}^{2},\quad
\Phi_{5}\,=\,\frac{19}{2^{3\mathstrut}\!\cdot\!3\!\cdot\!5^{2\mathstrut}}\,\sigma_{5}\,-\,\frac{1}{2^{2\mathstrut}\!\cdot\! 3^{2\mathstrut}}\,\sigma_{2}\sigma_{3},
\end{align*}
where
\[
\sigma_{k}\,=\,\sum\limits_{p\le y}\frac{1}{p^{\,k}},
\]
so that $\sigma = \sigma_{1}$ in this notation. The expressions for $\Phi_{n}$ become too complicated as $n$ grows.
The general algorithm for calculating  these coefficients is given in \cite[\S 1.2]{Korolev_2015}.
In particular, it is possible to show that $\Phi_{n}$ is a polynomial in variables $\sigma_{2}$, $\sigma_{3},\ldots, \sigma_{n}$.
For our purposes, we need only appropriate upper bounds for $|\Phi_{n}|$. Such bounds are given by the following
\vspace{0.5cm}

\textsc{Lemma 4.} \emph{The coefficients $\Phi_{n}$ satisfy the inequalities}

(a) $|\Phi_{n}|\le \frac{\displaystyle \sigma^{n}}{\displaystyle n!\mathstrut}$\;   \emph{for any} \; $n\ge 0$;

(b) $|\Phi_{n}|\le \frac{\displaystyle \sqrt{2\pi}}{\displaystyle n!\mathstrut}\,(\log\log{n}+1)^{n}$ \;  \emph{for} \;  $2\le n\le y$.

\vspace{0.5cm}

These are lemmas 1.11 and 1.12 from \cite{Korolev_2015}.

\vspace{0.5cm}

\textsc{Remark.} The quantity $\delta_{2k}$ defined in Lemma 1 can be expressed in terms of $\Phi_{n}$ as follows. Differentiating
both sides of (4.6) we have
\[
G^{(k)}(z)\,=\,\sum\limits_{n = 0}^{k}\binom{k}{n}\bigl(e^{\,\sigma z}\bigr)^{(n)}\Phi^{(k-n)}(z)\,=\,
\sum\limits_{n = 0}^{k}\binom{k}{n}\sigma^{n}e^{\,\sigma z}\Phi^{(k-n)}(z).
\]
Setting here $z = 0$ and using the relation $\Phi^{(k-n)}(0) = (k-n)!\Phi_{k-n}$, we obtain
\[
G^{(k)}(0)\,=\,k!\sum\limits_{n = 0}^{k}\frac{\Phi_{k-n}}{n!}\,\sigma^{n},\quad \delta_{2k}\,=\,
\frac{(2k)!}{k!}\,\frac{G^{(k)}(0)}{2^{2k}}\,=\,2^{-2k}\sum\limits_{n = 0}^{k}\frac{(2k)!}{n!}\,\Phi_{k-n}\,\sigma^{n}.
\]

\textsc{Lemma 5.} \emph{For any $u$ such that $|u|\le 1$ and for any integer $N\ge 1$ one has}
\[
G(-u^{2})\,=\,e^{\,-\sigma u^{2}}\biggl(\;\sum\limits_{n = 0}^{N}(-1)^{n}\Phi_{n}u^{2n}\,+\,O\bigl(R_{N}(u)\bigr)\biggr),
\]
\emph{where}
\[
R_{N}(u)\,=\,\frac{|u|^{2(N+1)}}{(N+1)!}\,\mu_{N}^{N+1},\quad \mu_{r}\,=\,\min{(\log\log{r}+1,\sigma)}.
\]

\textsc{Proof.} Since $G(z)$ is an entire function, then for any $u\in \mathbb{C}$ and for any $N\ge 1$ we get
\[
G(-u^{2})\,=\,e^{-\,\sigma u^{2}}\biggl(\;\sum\limits_{n = 0}^{N}(-1)^{n}\Phi_{n}u^{2n}\,+\,r_{N}(u)\biggr),\quad
r_{N}(u)\,=\,\sum\limits_{n>N}(-1)^{n}\Phi_{n}u^{2n}.
\]
Suppose now that $|u|\le 1$ and denote by $\nu$ the least integer satisfying the condition $\nu+1\ge y$. If $N\ge \nu$ then
in view of Lemma 4 (a) we have
\[
|r_{N}(u)|\,\le\,|u|^{2(N+1)}\sum\limits_{n = N+1}^{+\infty}a_{n},\quad a_{n}\,=\,\frac{\sigma^{n}}{n!}.
\]
Since
\[
\frac{a_{n+1}}{a_{n}}\,=\,\frac{\sigma}{N+1}\,<\,\frac{\sigma}{y}\,<\,\frac{1}{y}\,(\log\log{y}+1)\,<\,\frac{1}{8},
\]
then
\begin{equation*}\label{lab-2-6}
|r_{N}(u)|\,\le\,|u|^{2(N+1)}a_{N+1}\sum\limits_{k = 0}^{+\infty}\frac{1}{8^{k}}\,
=\,\frac{8}{7}\,\frac{|u|^{2(N+1)}}{(N+1)!}\,\sigma^{N+1}\,<\,\frac{8}{7}\,R_{N}(u).\eqno(4.7)
\end{equation*}
If $3\le N\le \nu-1$, then Lemma 4 (b) implies that
\[
|r_{N}(u)|\,\le\,\sqrt{2\pi}\sum\limits_{N<n\le \nu-1}b_{n}\,+\,|r_{\nu}(u)|,\quad b_{n}\,=\,\frac{1}{n!}\,(\log\log{n}+1)^{n}.
\]
To estimate the fraction $b_{n+1}/b_{n}$, we note that Lagrange's mean value theorem leads to the inequality
\[
\log\log{(n+1)}\,<\,\log\log{n}+(n\log{n})^{-1}
\]
for any $n\ge 2$. Therefore,
\begin{multline*}
\frac{b_{n+1}}{b_{n}}\,<\,\frac{\log\log{(n+1)}+1}{n+1}\,\biggl(\frac{\log\log{(n+1)}+1}{\log\log{n}+1}\biggr)^{\! n}\,<\\
<\,\frac{\log\log{(n+1)}+1}{n+1}\,\,\biggl(1\,+\,\frac{(n\log{n})^{-1}}{\log\log{n}+1}\biggr)^{\! n}\,<\,
\frac{\log\log{(n+1)}+1}{n+1}\,\exp{\biggl(\frac{(\log{n})^{-1}}{\log\log{n}+1}\biggr)}\,<\,\frac{4}{5}
\end{multline*}
for $n\ge 4$. Hence,
\[
\sum\limits_{N<n\le \nu-1}b_{n}\,<\,|u|^{2(N+1)}b_{N+1}\sum\limits_{k = 0}^{+\infty}\biggl(\frac{4}{5}\biggr)^{\! k}\,=\,
4|u|^{2(N+1)}b_{N+1}\,=\,5R_{N}(u).
\]
In view of (4.7) we have
\[
|r_{\nu}(u)|\,\le\,\frac{8}{7}\,\frac{|u|^{2(N+1)}}{(\nu +1)!}(\log\log{y}+1)^{\nu+1}\,=\,\frac{8}{7}\,R_{N}(u)\delta,
\]
where
\[
\delta\,=\,|u|^{2(\nu-N)}\,\frac{(N+1)!}{(\nu+1)!}\,\frac{(\log\log{y}+1)^{\nu+1}}{(\log\log{(N+1)}+1)^{N+1}}\,
\le\,\frac{b_{\nu+1}}{b_{N+1}}\,\le\,\biggl(\frac{4}{5}\biggr)^{\nu-N}\,\le\,\frac{4}{5}.
\]
Thus we obtain
\begin{equation*}
|r_{N}(u)|\,\le\,\biggl(5\sqrt{2\pi}\,+\,\frac{8}{7}\cdot\frac{4}{5}\biggr)R_{N}(u)\,<\,13.5R_{N}(u).\eqno(4.8)
\end{equation*}
If $N = 2$ then
\[
r_{2}(u)\,=\,-\Phi_{3}u^{6}\,+\,r_{3}(u).
\]
Using (4.8) together with the estimates
\[
|\Phi_{3}|\,<\,\frac{1}{9}\sum\limits_{p}\frac{1}{p^{3}}\,<\,\frac{1}{50},
\]
we find
\begin{alignat*}{3}
|r_{2}(u)|\,\le\,\frac{|u|^{6}}{50}\,+\,13.5R_{3}(u)\,\le\,\frac{|u|^{6}}{50}\,+\,13.5\cdot \frac{|u|^{8}}{4!}(\log\log{4}+1)^{4}\,\le\\
\le\,8.1\frac{|u|^{6}}{3!}(\log\log{3}+1)^{3}\,=\,8.1R_{2}(u).
\tag{4.9}
\end{alignat*}
Finally, if $N = 1$, then
\[
r_{1}(u)\,=\,\Phi_{2}u^{2}\,+\,r_{2}(u),
\]
so, using (4.9) together with the estimates
\[
|\Phi_{2}|\,<\,\frac{1}{4}\sum\limits_{p}\frac{1}{p^{2}}\,<\,\frac{1}{8},
\]
we obtain
\begin{multline*}
|r_{1}(u)|\,\le\,\frac{|u|^{4}}{8}\,+\,8.1R_{2}(u)\,\le\,\frac{|u|^{4}}{8}\,+\,8.1\cdot \frac{|u|^{6}}{3!}(\log\log{3}+1)^{3}\,\le\\
\le\,10\,\frac{|u|^{4}}{2!}(\log\log{2}+1)^{3}\,=\,10R_{1}(u).
\end{multline*}
The lemma is proved. $\Box$

\vspace{0.5cm}

\textsc{Lemma 6.} \emph{Suppose that the even integer $N\ge 2$ and real $z$ satisfy the inequality}
\[
\frac{\sigma e\pi^{2}(z^{2}+1)}{N}\,\le\,\frac{1}{8}.
\]
\emph{Then, for any $T$ and $H>1$, the integral}
\[
J(z)\,=\,J(z;T,H,y)\,=\,\int_{T}^{T+H}\exp{\bigl(2\pi iz V_{y}(t)\bigr)}\d t
\]
\emph{can be expressed in the form}
\[
J(z)\,=\,H G\bigl(-(\pi z)^{2}\bigr)\,+\,O\bigl(r(z)\bigr),\quad\textit{where}\quad
r(z)\,=\,|z|\biggl(\frac{H}{2^{N\mathstrut}\sqrt{N}}\,+\,y^{N/2}e^{e(\pi z)^{2}}\biggr).
\]

\textsc{Proof.} Using the notation of Lemma 1, we have
\begin{multline*}
J(z)\,=\,\int_{T}^{T+H}\biggl(\;\sum\limits_{n = 0}^{N-1}\frac{(2\pi iz)^{n}}{n!}\,V^{n}(t)\,+\,\theta\,\frac{(2\pi z)^{N}}{N!}\,V^{N}(t)\biggr)\d t\,=\\
=\,\sum\limits_{n = 0}^{N-1}\frac{(2\pi iz)^{n}}{n!}\,I_{n}\,+\,\theta_{1}\,\frac{(2\pi z)^{N}}{N!}\,I_{N}.
\end{multline*}
Denote by $R_{1}$ the contribution to $J(z)$ coming from the remainder term in (4.1). Since $I_{0} = H$, we have
\[
R_{1}\,\ll\,\sum\limits_{n = 1}^{N}\frac{(2\pi |z|)^{n}}{n!}\biggl(\frac{ny}{2}\biggr)^{\! n/2}\,\ll\,y^{N/2}\sum\limits_{n = 1}^{N}\lambda_{n}(z),\quad\text{where}\quad
\lambda_{n}(z)\,=\,\frac{(2\pi |z|)^{n}}{n!}\biggl(\frac{n}{2}\biggr)^{\! n/2}.
\]
By Stirling formula,
\[
\frac{s!\,s^{s}}{(2s)!}\,\ll\,\biggl(\frac{e}{4}\biggr)^{\! s}\,\ll\,\biggl(\frac{\sqrt{e}}{2}\biggr)^{\! 2s}
\]
for any $s\ge 1$. Hence, for odd $n = 2s+1$, $s\ge 0$, we obtain
\begin{multline*}
\lambda_{n}(z)\,=\,\frac{(2\pi|z|)^{2s+1}}{(2s+1)!}\,\bigl(s+0.5\bigr)^{s+0.5}\,\ll\,\frac{|z|}{\sqrt{s+1}}\,\frac{(2\pi|z|)^{2s}}{(2s)!}\,\bigl(s+0.5\bigr)^{s}\,\ll \\
\ll\,|z|\,\frac{(2\pi|z|)^{2s}}{s!}\,\frac{s!\,s^{s}}{(2s)!}\,\ll\,|z|\,\frac{(\pi\sqrt{e}\,|z|)^{2s}}{s!}.
\end{multline*}
Therefore, the contribution to $R_{1}$ coming from odd $n\le N$ does not exceed in order
\begin{equation*}
y^{N/2}|z|\sum\limits_{s = 0}^{+\infty}\frac{(\pi\sqrt{e}\,|z|)^{2s}}{s!}\,\ll\,|z|y^{N/2}e^{e(\pi z)^{2}}.\eqno(4.10)
\end{equation*}
Similarly, for even $n = 2s$, $s\ge 1$, and $|z|>1$ we have
\[
\lambda_{n}(z)\,=\,\frac{(2\pi|z|)^{2s}}{s!}\cdot\frac{s!\,s^{s}}{(2s)!}\,\ll\,\frac{(\pi\sqrt{e}\,|z|)^{2s}}{s!}\,\ll\,|z|\,\frac{(\pi\sqrt{e}|z|)^{2s}}{s!}.
\]
If $|z|\le 1$ then
\begin{multline*}
\lambda_{n}(z)\,\ll\,z^{2}\,\frac{(2\pi|z|)^{2(s-1)}}{(s-1)!}\,\frac{(s-1)!\,s^{s}}{(2s)!}\,\ll\,\frac{z^{2}}{s}\,\frac{(2\pi|z|)^{2(s-1)}}{(s-1)!}\,\cdot\frac{s!\,s^{s}}{(2s)!}\,\ll\\
\ll\,\frac{z^{2}}{s}\,\frac{(2\pi|z|)^{2(s-1)}}{(s-1)!}\,\biggl(\frac{\sqrt{e}}{2}\biggr)^{\! 2s}\,\ll\,\frac{z^{2}}{s}\,\frac{(\pi\sqrt{e}\,|z|)^{2(s-1)}}{(s-1)!}\,\ll\,
|z|\,\frac{(\pi\sqrt{e}\,|z|)^{2(s-1)}}{(s-1)!}.
\end{multline*}
Hence, the contribution to $R_{1}$ coming from even $n \le N$ is estimated as in (4.10). Therefore,
\[
R_{1}\,\ll\,|z|y^{N/2}e^{e(\pi z)^{2}}.
\]
Further, the sum of $\delta_{n}$ can be expressed as follows:
\begin{multline*}
\sum\limits_{0\le k\le N/2-1}\frac{(2\pi iz)^{2k}}{(2k)!}\cdot 2^{-2k}\,\frac{(2k)!}{k!}\,G^{(k)}(0)\,+
\,\theta\,\frac{(2\pi |z|)^{N}}{N!}\,2^{-N}\,\frac{N!}{(N/2)!}\,G^{(N/2)}(0)\,=\\
=\,\sum\limits_{k = 0}^{+\infty}\frac{(-(\pi z)^{2})^{k}}{k!}\,G^{(k)}(0)\,+\,R_{2}\,=\,G(-(\pi z)^{2})\,+\,R_{2},
\end{multline*}
where the term $R_{2}$ is estimated by Lemma 3 (c):
\[
R_{2}\,\ll\,\sum\limits_{k\ge N/2}\frac{(\pi|z|)^{2k}}{k!}\,G^{(k)}(0)\,\ll\,\sum\limits_{k\ge N/2}\frac{((\pi|z|)^{2}\sigma)^{k}}{k!}.
\]
If $k\ge N/2$ then Stirling formula implies the inequalities:
\[
\frac{((\pi|z|)^{2}\sigma)^{k}}{k!}\,\ll\,\frac{1}{\sqrt{k}}\biggl(\frac{\sigma e(\pi z)^{2}}{k}\biggr)^{\!k}\,\ll\,
\frac{1}{\sqrt{N}}\biggl(\frac{2\sigma e(\pi z)^{2}}{N}\biggr)^{\!k}.
\]
Hence,
\[
R_{2}\,\ll\,\frac{1}{\sqrt{N}}\biggl(\frac{2\sigma e(\pi z)^{2}}{N}\biggr)^{\!N/2}\;\sum\limits_{k = 0}^{+\infty}\biggl(\frac{2\sigma e(\pi z)^{2}}{N}\biggr)^{\!k}.
\]
Since
\[
\frac{2\sigma e(\pi z)^{2}}{N}\,<\,\frac{2\sigma e\pi^{2}(z^{2}+1)}{N}\,\le\,\frac{1}{4},
\]
then
\[
R_{2}\,\ll\,\frac{4^{-N/2}}{\sqrt{N}}\,\ll\,\frac{2^{-N}}{\sqrt{N\mathstrut}}\,\ll\,|z|\,\frac{2^{-N}}{\sqrt{N\mathstrut}}
\]
for $|z|>1$. In the case $|z|\le 1$ we have
\begin{multline*}
R_{2}\,\ll\,\frac{1}{\sqrt{N}}\biggl(\frac{2\sigma e(\pi z)^{2}}{N}\biggr)^{\!N/2}\,\ll\,\frac{|z|^{N}}{\sqrt{N\mathstrut}}
\biggl(\frac{2\sigma e\pi^{2}}{N}\biggr)^{\!N/2}\,\ll\\
\ll\,\frac{|z|^{N}}{\sqrt{N\mathstrut}}
\biggl(\frac{2\sigma e\pi^{2}(z^{2}+1)}{N}\biggr)^{\!N/2}\,\ll\,|z|^{N}\,\frac{4^{-N/2}}{\sqrt{N\mathstrut}}\,\ll\,|z|^{N}\,\frac{2^{-N}}{\sqrt{N\mathstrut}}\,\ll\,
|z|\,\frac{2^{-N}}{\sqrt{N\mathstrut}}.
\end{multline*}
Now the lemma is proved. $\Box$

\vspace{0.5cm}

\textsc{Lemma 7.} \emph{For any $c\ge 1$ the following inequality holds:}
\[
\int_{1}^{+\infty}\frac{e^{-cv^{2}}}{v}\d v\,\le\,\frac{e^{-c}}{2c}.
\]

\textsc{Proof.} Taking $v = \sqrt{u/c}$ in the integral, we get
\[
\int_{1}^{+\infty}\frac{e^{-cv^{2}}}{v}\d v\,=\,\frac{1}{2}\int_{c}^{+\infty}\frac{e^{-u}}{u}\d u\,\le\,\frac{1}{2c}\int_{c}^{+\infty}e^{-u}\d u\,=\,\frac{e^{-c}}{2c}.
\]
Lemma 7 is proved. $\Box$

 \medskip

\textsc{Lemma 8.} \emph{If $H_n(x)$ is the Hermite polynomial, then for any $n\ge 0$ the following relations hold}:
\begin{align*}
\text{(a)} & \quad \int_{-\infty}^{+\infty}x^{2n}e^{-\,x^{2}/2}e^{ixy}\d x\,=\,\frac{(-1)^{n}}{2^{2n\mathstrut}}\sqrt{\frac{\pi}{2}}\,e^{-\,y^{2}/2}H_{2n}(y),\\
\text{(b)} & \quad \int_{0}^{+\infty}H_{2n}(x)e^{-\,x^{2}/2}\cos{(xy)}\d x\,=\,(-1)^{n}\sqrt{\frac{\pi}{2}}\;y^{2n}e^{-\,y^{2}/2},\\
\text{(c)} & \quad \int_{0}^{+\infty}e^{-\,x^{2}/2}H_{n}^{2}(x)\d x\,=\,\sqrt{\frac{\pi}{2}}\,n!.
\end{align*}

These relations are well-known (see, for example, \cite[\S\S 1,2, Ch. V]{Suetin_2005}).

\vspace{0.5cm}

\textsc{Lemma 9.} \emph{For any $n\ge 0$, the following inequality holds}:
\[
\int_{0}^{+\infty}e^{-\,x^{2}/2}|H_{2n}(x)|\d x\,\le\,\sqrt{\frac{\pi}{2}}\,\sqrt{(2n)!}\,.
\]

\textsc{Proof.} Denoting the integral above by $\kappa_{n}$  and using Cauchy's inequality together with Lemma 8 (c), we get:
\[
\kappa_{n}^{2}\,\le\,\biggl(\;\int_{0}^{+\infty}e^{-\,x^{2}/2}\d x\biggr)\int_{0}^{+\infty}e^{-\,x^{2}/2}H_{2n}^{2}(x)\d x
\,\le\,\sqrt{\frac{\pi}{2}}\,\sqrt{\frac{\pi}{2}}\,(2n)!\,=\,
\frac{\pi}{2}\,(2n)!.
\]
Lemma 9 is proved. $\Box$

\vspace{0.5cm}

\textsc{Lemma 10.} \emph{Suppose that $n\ge  0$ is integer, $\omega>1$, $\lambda>0$, and $K(x)$ is defined by} (3.2).
\emph{Then, for any real $\xi$, the integral}
\[
j_{n}\,=\,j_{n}(\omega,\lambda,\xi)\,=\,\int_{-\infty}^{+\infty}|H_{2n}(u)|e^{-\,u^{2}/2}\,K\bigl(\omega(\xi+\lambda u)\bigr)\d u
\]
\emph{satisfies the inequality}:
\[
j_{n}\,\le\,\frac{(n+1)^{1/4}}{\omega\lambda\sqrt{2}}\,\frac{(2n)!}{n!}\,\biggl(\frac{3}{2}\biggr)^{\!n}.
\]
\textsc{Proof.} By using (3.3)  we obtain
\[
j_{n}\,\le\,\int_{-\infty}^{+\infty}\sum\limits_{r = 0}^{n}\frac{(2n)!}{(2r)!}\,\frac{(2u)^{2r}}{(n-r)!}
\,e^{-\,u^{2}/2}K\bigl(\omega(\xi+\lambda u)\bigr)\d u\,=\,
(2n)!\sum\limits_{r = 0}^{n}\frac{2^{2r}}{(n-r)!}\,\frac{k_{r}}{(2r)!},
\]
where
\[
k_{r}\,=\,\int_{-\infty}^{+\infty}u^{2r}e^{-\,u^{2}/2}\,K\bigl(\omega(\xi+\lambda u)\bigr)\d u.
\]
We shall use now the well-known identity
\begin{equation*}
K(x)\,=\,2\int_{0}^{1}(1-t)\cos{(2\pi xt)}\d t.\eqno(4.11)
\end{equation*}
It follows that
\begin{multline*}
k_{r}\,=\,\RRe\int_{-\infty}^{+\infty}u^{2r}e^{-\,u^{2}/2}\,\frac{2}{\omega}
\int_{0}^{\,\omega}\biggl(1-\frac{y}{\omega}\biggr)e^{2\pi iy(\xi+\lambda u)}\d y\,\d u\,=\\
=\,\frac{2}{\omega}\RRe\int_{0}^{\,\omega}\biggl(1-\frac{y}{\omega}\biggr)e^{2\pi iy\xi}
\int_{-\infty}^{+\infty}u^{2r}e^{-\,u^{2}/2+2\pi i\lambda yu}\d u\d y.
\end{multline*}
By Lemma 8 (a), we have
\[
k_{r}\,=\,\frac{2}{\omega}\RRe\int_{0}^{\,\omega}\biggl(1-\frac{y}{\omega}\biggr)
e^{2\pi iy\xi}\cdot \frac{(-1)^{r}}{2^{2r}}\sqrt{\frac{\pi}{2}}\,e^{-2(\pi\lambda y)^{2}}H_{2r}(2\pi\lambda y)\d y,
\]
and hence, by Lemma 9,
\begin{multline*}
|k_{r}|\,\le\,\frac{2}{\omega 2^{2r}}\sqrt{\frac{\pi}{2}}\int_{0}^{+\infty}e^{-2(\pi\lambda y)^{2}}\bigl|H_{2r}(2\pi\lambda y)\bigr|dy\,=\,
\frac{(\omega\lambda)^{-1}}{2^{2r}\sqrt{2\pi}}\int_{0}^{+\infty}e^{-\,x^{2}/2}\bigl|H_{2r}(x)\bigr|dx\,\le\\
\le\,\frac{(\omega\lambda)^{-1}}{2^{2r}\sqrt{2\pi}}\,\sqrt{\frac{\pi}{2}}\sqrt{(2r)!}\,=\,
\frac{(\omega\lambda)^{-1}}{2^{2r+1}}\,\sqrt{(2r)!}\,.
\end{multline*}
Therefore,
\begin{multline*}
j_{n}\,\le\,(2n)!\sum\limits_{r = 0}^{n}\frac{2^{2r}}{(n-r)!}\,\frac{(\omega\lambda)^{-1}}{2^{2r+1}}\,\sqrt{(2r)!}\,=\,(2\omega\lambda)^{-1}\,\frac{(2n)!}{n!}
\sum\limits_{r = 0}^{n}\frac{n!}{r!(n-r)!}\binom{2r}{r}^{-1/2}\,\le\\
\le\,(2\omega\lambda)^{-1}\,\frac{(2n)!}{n!}\sum\limits_{r = 0}^{n}\binom{n}{r}\cdot\frac{(r+1)^{1/4}}{2^{2r+0.5}}\,\le\,
\frac{(n+1)^{1/4}}{\omega\lambda\sqrt{2}}\,\frac{(2n)!}{n!}\sum\limits_{r = 0}^{n}2^{-r}\binom{n}{r}\,=\\
=\,\frac{(n+1)^{1/4}}{\omega\lambda\sqrt{2}}\,\frac{(2n)!}{n!}\,\biggl(\frac{3}{2}\biggr)^{\!n}.
\end{multline*}
The lemma is proved.

\vspace{0.5cm}

We define the function $g(t)$ as follows:
\begin{equation*}
g(t)\,=\,
\begin{cases}
\displaystyle 2\bigl(\pi^{-1}\text{sgn}(t)\,+\,(1-|t|)\cot(\pi t)\bigr), & \text{for}\;|t|\le 1,\\
0, & \text{otherwise}.
\end{cases}
\end{equation*}
One can check that $g(t)$ is an odd, unbounded function such that
\begin{align*}
& g(t)\,=\,\frac{1}{\pi t}\,-\,\frac{2\pi}{3}\,t\,+\,\frac{2\pi}{3}\,t^{2}\,+\,O\bigl(t^{3}\bigr), \quad \text{as}\;t\to 0+,\\
& g(t)\,=\,\frac{\pi}{3}(1-t)^{2}\,+\,O\bigl((1-t)^{4}\bigr), \quad \text{as}\;t\to 1-0.
\end{align*}
Further, for any $\omega>0$ we set
\[
F_{\omega}(u)\,=\,\frac{1}{\omega}\int_{0}^{\,\omega}g\biggl(\frac{t}{\omega}\biggr)\sin{(2\pi ut)}\d t.
\]
Then we have
\vspace{0.5cm}

\textsc{Lemma 11.} \emph{Let $K(x)$ be as in} (3.2). \emph{For any real $u$ and any $\omega>0$ one has}
\[
\text{sgn}(u)\,=\,F_{\omega}(u)\,+\,\theta K(\omega u).
\]

This is lemma 4.1 from \cite{Tsang_1984}. Note that our definition of $F_{\omega}(u)$
differs slightly from the definition introduced in \cite{Tsang_1984},
but it serves a similar purpose as in \cite{Tsang_1984}.
The fact that $|\theta| \le 1$ follows from the inequalities at the bottom of p. 28 of \cite{Tsang_1984}.
\vspace{0.5cm}

\textsc{Remark.} The functions that approximate the sign-function $\text{sgn}(u)$,
or the characteristic function $\chi_{E}(u)$ of any segment $E\subset \mathbb{R}$,
were discovered independently by A.~Selberg and A.~Beurling. They are of a great importance in approximation theory.
For an extensive account, see the paper of J.~Vaaler \cite{Vaaler_1985}.

\vspace{0.5cm}
\textsc{Lemma 12.} \emph{Suppose that $0<\varepsilon<10^{-3}$ is an arbitrary small fixed constant,
$T\ge T_{0}(\varepsilon)$, $x=T^{0.1\varepsilon}$, $T^{c+\varepsilon_{1}}\le H\le T^{c+\varepsilon}$,
where $c = \tfrac{27}{82}$, $\varepsilon_{1} =0.9\varepsilon$, and let $1\le m\le c_{1}\log{x}$, where $c_{1}$ is a sufficiently small absolute constant. Further,
let $y = x^{1/(8m+3)}$. Then the following inequality holds:}
\[
\int_{T}^{T+H}\bigl|\pi S(t)\,+\,V_{y}(t)\bigr|^{2m}\d t\,\le\,\frac{c_{0}H}{56\varepsilon}\,(\pi c_{0}\varepsilon^{-1})^{2m},
\]
\emph{where} $c_{0} = 2880$.

\vspace{0.5cm}

This is lemma 3.13 from \cite{Korolev_2015}. The lemma provides one way of showing that $S(t)$ is well approximated
by $-\frac{1}{\pi}V_{y}(t)$ over the short interval $[T, T+H]$, which is of crucial importance for the proofs of our results. It also shows how the constant $c = \tfrac{27}{82}$ appears in our results.

\vspace{0.5cm}

\textbf{5. Proof of Theorem 3}

\vspace{0.5cm}
As stated in Section 3, Theorem 3 is the fundamental result which will enable us to deduce Theorem 2. Thus we
start with the proof of this result.

\medskip
 We shall follow the proof of Theorem 6.1 from \cite{Tsang_1984}, with appropriate changes.
 Let $\omega$ and the even integer $N$ satisfy the conditions
\begin{equation*}
1\,\le\,\omega\,\le\,\frac{1}{8\pi\sqrt{e}}\sqrt{\frac{\log{H}}{\sigma\log{y}}},\quad 1\,\le\,
\omega\,\le\,\frac{1}{4}\,\exp{\bigl(0.5(\log{y})^{1/3}\bigr)},\eqno(5.1)
\end{equation*}
\begin{equation*}
2\,\le\,N\,\le\,\frac{\log{H}}{4\log{y}},\quad \frac{\sigma\omega^{2}}{N}\,\le\,\frac{1}{16\pi^{2}e}.\eqno(5.2)
\end{equation*}
Note that such pairs $\omega, N$ exist. Indeed, from (2.4) we have $y<H$ and
\[
\frac{\log{H}}{\log{y}}\,>\,14\,000\log\log{H}\,>\,2(16\pi)^{2}e\log\log{y}\,>\,(16\pi)^{2}e(\log\log{y}+1)\,>\,(16\pi)^{2}e\sigma.
\]
Therefore,
\[
\frac{1}{8\pi\sqrt{e}}\sqrt{\frac{\log{H}}{\sigma\log{y}}}\,>\,1.
\]
Next, since
\[
\frac{\log{H}}{8\log{y}}\,>\,\frac{1}{8}\cdot 14\,000\log\log{H}\,>\,1\,750\log\log{H},
\]
there are at least $825\log\log{H}$ even integers between
\[
\frac{\log{H}}{8\log{y}}\quad\text{and}\quad \frac{\log{H}}{4\log{y}}.
\]
Now let $\omega$ be any number from the segment
\[
1\,\le\,\omega\,\le\,\min{\biggl\{\frac{1}{8\pi\sqrt{e}}\sqrt{\frac{\log{H}}{\sigma\log{y}}},\frac{1}{4}\,\exp{\bigl(0.5(\log{y})^{1/3}\bigr)}\biggr\}},
\]
and $N$ be any even integer from the segment
\[
\frac{\log{H}}{8\log{y}}\,\le\,N\,\le\,\frac{\log{H}}{4\log{y}}.
\]
Then we have
\[
N\,\ge\,\frac{\log{H}}{8\log{y}}\,=\,\frac{\log{H}}{\sigma\log{y}}\,\frac{\sigma}{8}\,\frac{(16\pi)^{2}e}{(16\pi)^{2}e}\,=\,32\pi^{2}e\sigma\cdot
\frac{1}{(16\pi)^{2}e}\,\frac{\log{H}}{\sigma\log{y}}\,\ge\,32\pi^{2}e\sigma\omega^{2}
\]
and therefore
\[
\frac{\sigma\omega^{2}}{N}\,\le\,\frac{1}{32\pi^{2}e}.
\]
Further, setting $V(t) = V_{y}(t)$ for brevity and applying Lemma 11, we get
\[
\int_{T}^{T+H}\text{sgn}(V(t)-\alpha)\d t\,=\,I_{1}\,+\,\theta R_{1},
\]
where
\[
I_{1}\,=\,\int_{T}^{T+H}F_{\omega}\bigl(V(t)-\alpha\bigr)\d t,\quad R_{1}\,=\,\int_{T}^{T+H}K\bigl(\omega(V(t)-\alpha)\bigr)\d t.
\]
First we estimate $R_{1}$. By the identity (4.11) we have
\begin{alignat*}{4}
K(\omega x)\,=\,2\int_{0}^{1}(1-y)\cos{(2\pi x\omega y)}\d y
\,=\,\frac{2}{\omega}\int_{0}^{\,\omega}\biggl(1-\frac{v}{\omega}\biggr)\cos{(2\pi xv)}\d v\,=\\
=\,\frac{2}{\omega}\,\RRe\int_{0}^{\,\omega}\biggl(1-\frac{v}{\omega}\biggr)e^{2\pi ixv}\d v.
\tag{5.3}
\end{alignat*}
Taking $x = V(t)-\alpha$ in (5.3), we express $R_{1}$ as follows:
\begin{multline*}
R_{1}\,=\,\frac{2}{\omega}\,\RRe\int_{T}^{T+H}\int_{0}^{\,\omega}\biggl(1-\frac{v}{\omega}\biggr)e^{2\pi iv(V(t)-\alpha)}\d v\d t\,=\\
=\,\frac{2}{\omega}\,\RRe\int_{0}^{\,\omega}\biggl(1-\frac{v}{\omega}\biggr)e^{-2\pi iv\alpha}\int_{T}^{T+H}e^{2\pi ivV(t)}\d t\d v\,=\\
=\,\frac{2}{\omega}\,\RRe\int_{0}^{\,\omega}\biggl(1-\frac{v}{\omega}\biggr)e^{-2\pi iv\alpha}J(v)\d v,
\end{multline*}
where $J(v)$ is as in Lemma 6. By the conditions (5.1), and (5.2), for any $v$, $0\le v\le \omega$, we have
\[
0\,<\,\frac{\sigma(v^{2}+1)}{N}\,\le\,\frac{\sigma(\omega^{2}+1)}{N}\,\le\,
\frac{2\sigma\omega^{2}}{N}\,\le\,\frac{1}{16\pi^{2}e},
\]
and hence
\[
\frac{\sigma e\pi^{2}(v^{2}+1)}{N}\,\le\,\frac{1}{16}\,<\,\frac{1}{8}.
\]
Thus Lemma 6 yields:
\begin{equation*}
J(v)\,=\,HG\bigl(-(\pi v)^{2}\bigr)\,+\,O\biggl(|v|\biggl(\,H\,\frac{2^{-N}}{\sqrt{N}}\,+\,y^{N/2}e^{e(\pi v)^{2}}\biggr)\biggr).\eqno(5.4)
\end{equation*}
It is easy to check that
\begin{multline*}
e(\pi v)^{2}\,\le\,e(\pi\omega)^{2}\,\le\,e\pi^{2}\,\frac{1}{(8\pi)^{2}e}\,\frac{\log{H}}{\sigma\log{y}}
\,=\,\frac{\log{H}}{64\sigma\log{y}}\,<\,\frac{1}{64}\,\log{H},\\
y^{N/2}\,\le\,\exp{\biggl(\frac{1}{2}\,\frac{\log{H}}{4\log{y}}\,\log{y}\biggr)}\,=\,H^{1/8}.
\end{multline*}
Therefore we obtain
\[
y^{N/2}e^{e(\pi v)^{2}}\,\le\,H^{1/8+1/64}\,<\,H^{1/7}.
\]
At the same time,
\[
\frac{2^{-N}}{\sqrt{N}}\,\ge\,\frac{1}{\sqrt{N}}\,\exp{\biggl(-\,\frac{1}{4}\,\frac{\log{H}}{\log{y}}\,\log{2}\biggr)}
\,>\,\exp{\biggl(-\,\frac{1}{6}\,\log{H}\biggr)}\,=\,H^{-1/6}.
\]
Hence, if $0\le v\le\omega$ then
\begin{equation*}
y^{N/2}e^{e(\pi v)^{2}}\,\le\,H\cdot H^{-6/7}\,<\,H\cdot H^{-1/6}\,<\,\frac{H}{2^{N\mathstrut}\sqrt{N}}.\eqno(5.5)
\end{equation*}
Therefore, the $O$-term in (5.4) contributes to $R_{1}$ at most
\[
\frac{1}{\omega}\int_{0}^{\,\omega}v\,H\,\frac{2^{-N}}{\sqrt{N}}\d v\,\ll\,\frac{\omega H}{2^{N\mathstrut}\sqrt{N}}.
\]
Thus we have
\begin{multline*}
R_{1}\,=\,\frac{2H}{\omega}\,\RRe \int_{0}^{\,\omega}\biggl(1-\frac{v}{\omega}\biggr)G\bigl(-(\pi v)^{2}\bigr)e^{-2\pi i\alpha v}\d v
\,+\,O\biggl(\omega H\,\frac{2^{-N}}{\sqrt{N}}\biggr)\,=\\
=\,\frac{2H}{\pi\omega}\,\RRe \int_{0}^{\,\pi\omega}\biggl(1-\frac{v}{\pi\omega}\biggr)G\bigl(-v^{2}\bigr)e^{-2i\alpha v}\d v
\,+\,O\biggl(\omega H\,\frac{2^{-N}}{\sqrt{N}}\biggr).
\end{multline*}
Lemmas 3 (b) and 5 imply the estimate $G(-v^{2})\ll e^{-\sigma v^{2}/2}$ for any $v$, $0\le v\le \pi \omega$. Thus,
\[
R_{1}\,\ll\,\frac{H}{\omega}\int_{0}^{\,\pi\omega}e^{-\sigma v^{2}/2}\d v\,+\,\frac{\omega H}{2^{N\mathstrut}\sqrt{N}}\,\ll\,\frac{H}{\omega\sqrt{\sigma}}\,+\,
\frac{\omega H}{2^{N\mathstrut}\sqrt{N}}.
\]
In view of (5.2), one has $N>\omega^{2}\sigma$, hence,
\begin{equation*}
2^{N}\sqrt{N}\,>\,2^{\omega^{2}\sigma}\omega\sqrt{\sigma}\,>\,\omega^{2}\sqrt{\sigma}\quad\text{and}\quad\frac{\omega}{2^{N\mathstrut}\sqrt{N}}
\,<\,\frac{1}{\omega\sqrt{\sigma}}.\eqno(5.6)
\end{equation*}
Finally, we conclude that
\[
R_{1}\,\ll\,\frac{H}{\omega\sqrt{\sigma}}.
\]
Now we calculate the integral $I_{1}$. First, we have
\begin{align*}
I_{1}\,&=\,\frac{1}{\omega}\int_{T}^{T+H}\int_{0}^{\,\omega}g\biggl(\frac{x}{\omega}\biggr)\sin{\bigl\{2\pi x(V(t)-\alpha)\bigr\}}\d x\d t\\&
=\frac{1}{\omega}\int_{0}^{\,\omega}g\biggl(\frac{x}{\omega}\biggr)\IIm \int_{T}^{T+H}e^{2\pi ix(V(t)-\alpha)}\d t\d x,\\&
=\frac{1}{\omega}\int_{0}^{\,\omega}g\biggl(\frac{x}{\omega}\biggr)\IIm\biggl(e^{-2\pi i\alpha x}\int_{T}^{T+H}e^{2\pi ixV(t)}\d t\biggr)\d x\\&
=\,\frac{1}{\omega}\int_{0}^{\,\omega}g\biggl(\frac{x}{\omega}\biggr)\IIm\bigl(J(x)e^{-2\pi i\alpha x}\bigr)\d x.
\end{align*}
The application of Lemma 6 gives: $I_{1} = I_{2}+O(R_{2})$, where
\begin{align*}
I_{2}\,&=\,\frac{1}{\omega}\int_{0}^{\,\omega}g\biggl(\frac{x}{\omega}\biggr)\IIm\bigl(HG\bigl(-(\pi x)^{2}\bigr)e^{-2\pi i\alpha x}\bigr)\d x\\&
=\,-\frac{H}{\omega}\int_{0}^{\,\omega}g\biggl(\frac{x}{\omega}\biggr)\sin{(2\pi\alpha x)}G\bigl(-(\pi x)^{2}\bigr)\d x\\&\,=\,
-\frac{H}{\pi\omega}\int_{0}^{\,\pi\omega}g\biggl(\frac{x}{\pi\omega}\biggr)\sin{(2\alpha x)}G\bigl(-x^{2}\bigr)\d x,\\&
R_{2}\,=\,\frac{1}{\omega}\int_{0}^{\,\omega}g\biggl(\frac{x}{\omega}\biggr)x\biggl(\,\frac{\omega 2^{-N}}{\sqrt{N}}\,+\,y^{N/2}e^{e(\pi x)^{2}}\biggr)\d x.
\end{align*}
In view of (5.5) and (5.6), the term $R_{2}$ is estimated as follows:
\[
R_{2}\,\ll\,\frac{H}{\omega}\,\frac{2^{-N}}{\sqrt{N}}\int_{0}^{\,\omega}g\biggl(\frac{x}{\omega}\biggr)x\d x\,\ll\,
\frac{H}{\omega}\,\frac{2^{-N}}{\sqrt{N}}\int_{0}^{\,\omega}\frac{\omega}{x}\,x\d x
\,\ll\,\frac{\omega H}{2^{N\mathstrut}\sqrt{N}}\,\ll\,\frac{H}{\omega\sqrt{\sigma}}.
\]
Next, by Lemmas 3 (b) and 7, the integral over the segment $1\le x\le \pi \omega$ in the expression for $I_{2}$ is bounded by
\[
\frac{H}{\omega}\int_{1}^{\,\pi\omega}g\biggl(\frac{x}{\pi\omega}\biggr)e^{-\,\sigma x^{2}/2}\d x
\,\ll\,\frac{H}{\omega}\int_{1}^{\,\pi\omega}\frac{\omega}{x}\,e^{-\,\sigma x^{2}/2}\d x\,\ll\,H
\int_{1}^{+\infty}e^{-\,\sigma x^{2}/2}\,\frac{\d x}{x}\,\ll\,\frac{H}{\sigma}\,e^{-\,\sigma/2}.
\]
For $0\le x\le 1$, we expand $G(-x^{2})$ into Taylor series by Lemma 5. Thus we obtain
\begin{multline*}
I_{2}\,=\\
=\,-\frac{H}{\pi\omega}\int_{0}^{1}g\biggl(\frac{x}{\pi\omega}\biggr)e^{-\,\sigma x^{2}}\sin{(2\alpha x)}\biggl(\;\sum\limits_{n=0}^{\nu}(-1)^{n}\Phi_{n}x^{2n}
\,+\,O\bigl(x^{2(\nu+1)}D_{\nu}\bigr)\biggr)\d x\,+\,O\biggl(\frac{H}{\sigma}\,e^{-\,\sigma/2}\biggr),
\end{multline*}
where
\[
D_{\nu}\,=\,\frac{(\log\log{(\nu+2)}+1)^{\!\nu+1}}{(\nu+1)!}.
\]
The contribution to the integral coming from the $O$-term in the integrand is estimated as
\begin{multline*}
\ll\,\frac{H}{\omega}\,D_{\nu}\int_{0}^{1}g\biggl(\frac{x}{\pi\omega}\biggr)x^{2\nu+2}e^{-\,\sigma x^{2}}\d x\,\ll\,
\frac{H}{\omega}\,D_{\nu}\int_{0}^{1}\frac{\omega}{x}\,x^{2\nu+2}e^{-\,\sigma x^{2}}\d x\,\ll\\
\ll\,HD_{\nu}\int_{0}^{1}x^{2\nu+1}e^{-\,\sigma x^{2}}\d x
\,\ll\,\frac{HD_{\nu}}{\sigma^{\nu+1}}\int_{0}^{+\infty}u^{\nu}e^{-u}\d u\,\ll\,\nu!\,\frac{HD_{\nu}}{\sigma^{\nu+1}}\,\ll\\
\ll\,\frac{H}{\nu+1}\biggl(\frac{\log\log{(\nu+2)}+1}{\sigma}\biggr)^{\!\nu+1}
\end{multline*}
This means that $I_{2} = I_{3} + O(R_{3})$, where
\begin{equation*}
I_{3}\,=\,-\frac{H}{\pi\omega}\int_{0}^{1}g\biggl(\frac{x}{\pi\omega}\biggr)
e^{-\,\sigma x^{2}}\sin{(2\alpha x)}\biggl(\;\sum\limits_{n=0}^{\nu}(-1)^{n}\Phi_{n}x^{2n}\biggr)\d x, \eqno(5.7)
\end{equation*}
\[
R_{3}\,=\,\frac{H}{\sigma}\,e^{-\sigma/2}\,+\,\frac{H}{\nu+1}\biggl(\frac{\log\log{(\nu+2)}+1}{\sigma}\biggr)^{\!\nu+1}.
\]
Now we replace the limits in the integral (5.7) by $0\le x\le \pi\omega$. First we note that if $2\le n\le\nu$, then
\begin{alignat*}{5}
&\frac{H}{\omega}\int_{1}^{\pi\omega}g\biggl(\frac{x}{\pi\omega}\biggr)x^{2n}e^{-\,\sigma x^{2}}
\d x&\,\ll\,H\int_{1}^{+\infty}x^{2n-1}e^{-\,\sigma x^{2}}\d x\,\ll\\&
\ll\,He^{-\,\sigma/2}\int_{1}^{+\infty}x^{2n-1}e^{-\,\sigma x^{2}/2}\d x\,\ll
\,He^{-\,\sigma/2}\biggl(\frac{2}{\sigma}\biggr)^{\!n}(n-1)!\,.
\tag{5.8}
\end{alignat*}

Next, by Lemma 7 for $n=0$, we get
\begin{equation*}
\frac{H}{\omega}\int_{1}^{\pi\omega}g\biggl(\frac{x}{\pi\omega}\biggr)e^{-\,\sigma x^{2}}\d x
\,\ll\,H\int_{1}^{\,\pi\omega}\frac{e^{-\,\sigma x^{2}}}{x}\d x\,\ll\,
\frac{H}{\sigma}\,e^{-\sigma}.\eqno(5.9)
\end{equation*}
Using (5.8) and (5.9), we find that the above change of limits of integration contributes to (5.7) at most
\[
\frac{H}{\omega}\sum\limits_{n=0}^{\nu}|\Phi_{n}|\int_{1}^{\,\pi\omega}g\biggl(\frac{x}{\pi\omega}\biggr)x^{2n}e^{-\,\sigma x^{2}}\d x\,\ll\,
H\biggl(\,\frac{e^{-\sigma}}{\sigma}\,+\,\sum\limits_{n=2}^{\nu}|\Phi_{n}|e^{-\sigma/2}\biggl(\frac{2}{\sigma}\biggr)^{\!n}(n-1)!\biggr).
\]
By Lemma 3 (b), the last sum is $\ll e^{-\sigma/2}\Sigma_{\nu}$, where
\[
\Sigma_{\nu}\,\ll\,\sum\limits_{n = 2}^{\nu}\frac{1}{n}\biggl(\frac{2}{\sigma}\,(\log\log{n}+1)\biggr)^{\!n}.
\]
Set $\nu_{0} = \bigl[\exp{\bigl((\log{y})^{1/4}\bigr)}\bigr]$, $\nu_{1} = [\log{y}]$. If $2\le\nu\le \nu_{1}$, then
\begin{multline*}
\Sigma_{\nu}\,\ll\,\frac{1}{\sigma^{2}}\,+\,\frac{1}{\sigma^{3}}\,+
\,\sum\limits_{n = 4}^{\nu}\frac{1}{n}\biggl(\frac{2}{\sigma}\,(\log\log{\nu_{1}}+1)\biggr)^{\!n}\,\ll\,
\frac{1}{\sigma^{2}}\,+\,\sum\limits_{n = 4}^{\nu}\frac{1}{n}\biggl(\frac{2(\log{\sigma}+1)}{\sigma}\biggr)^{\!n}\,\ll\\
\ll\,\frac{1}{\sigma^{2}}\,+\,\biggl(\frac{2(\log{\sigma}+1)}{\sigma}\biggr)^{\!4}\sum\limits_{n = 4}^{\nu}\frac{1}{n}\,
\ll\,\frac{1}{\sigma^{2}}\,+\,\biggl(\frac{\log{\sigma}}{\sigma}\biggr)^{\!4}\cdot\sigma\,\ll\,\frac{1}{\sigma^{2}}.
\end{multline*}
If $\nu_{1}<\nu\le \nu_{0}$, then we estimate the sum $\Sigma_{\nu}$ as follows:
\[
\Sigma_{\nu}\,\ll\,\Sigma_{\nu_{1}}\,+\,\sum\limits_{n = \nu_{1}+1}^{\nu_{0}}\frac{1}{n}\biggl(\frac{2}{\sigma}\,(\log\log{\nu_{1}}+1)\biggr)^{\!n}\,\ll\,
\frac{1}{\sigma^{2}}\,+\,\sum\limits_{n = \nu_{1}+1}^{\nu_{0}}\frac{1}{n}\biggl(\frac{2}{\sigma}\,(\log\log{\nu_{0}}+1)\biggr)^{\!n}.
\]
Since
\[
\frac{2}{\sigma}\,(\log\log{\nu_{0}}+1)\,\le\,\frac{2}{\sigma}\,\biggl(\frac{1}{4}\log\log{y}+1\biggr)\,\le\,\frac{2}{\sigma}\biggl(\frac{\sigma}{4}+1\biggr)\,=\,
\frac{1}{2}+\frac{2}{\sigma}\,<\,\frac{2}{3},
\]
then
$$
\Sigma_{\nu} \ll \frac{1}{\sigma^{2}}\,+\,\sum\limits_{n = \nu_{1}+1}^{\nu_{0}}\frac{1}{n}\biggl(\frac{2}{3}\biggr)^{\!n}\ll
\frac{1}{\sigma^{2}}\,+\,\frac{1}{\nu_{1}}\biggl(\frac{2}{3}\biggr)^{\!\nu_{1}}\,
\ll\,\frac{1}{\sigma^{2}}\,+\,\frac{y^{-\log{(3/2)}}}{\log{y}}\,\ll\,\frac{1}{\sigma^{2}}\,+\,y^{-2/5}\,\ll\,\frac{1}{\sigma^{2}}.
$$
Thus, for any $\nu$, $0\le \nu\le \nu_{0}$, we have
\begin{alignat*}{6}
I_{3} &=-\frac{H}{\pi\omega}\int_{0}^{\,\pi\omega}g\biggl(\frac{x}{\pi\omega}\biggr)e^{-\sigma x^{2}}\sin{(2\alpha x)}
\biggl(\;\sum\limits_{n=0}^{\nu}(-1)^{n}\Phi_{n}x^{2n}\biggr)\d x\,+\,O\biggl(\frac{He^{-\,\sigma/2}}{\sigma^{2}}\biggr)\\&
=\,\frac{H}{\pi}\sum\limits_{n=0}^{\nu}(-1)^{n}\Phi_{n}\,j(n)\,+\,O\biggl(\frac{He^{-\,\sigma/2}}{\sigma^{2}}\biggr),
\tag{5.10}
\end{alignat*}
where we set
\begin{alignat*}{7}
j(n)\,=\,-\frac{1}{\omega}\int_{0}^{\,\pi\omega}g\biggl(\frac{x}{\pi\omega}\biggr)
\sin{(2\alpha x)}x^{2n}e^{-\sigma x^{2}}\d x.\tag{5.11}
\end{alignat*}
By Lemma 8 (b), for any $y>0$ we have
\[
y^{2n}e^{-y^{2}/2}\,=\,(-1)^{n}\sqrt{\frac{2}{\pi}}\int_{0}^{+\infty}H_{2n}(v)e^{-\,v^{2}/2}\cos{(vy)}\d v.
\]
Taking $y = x\sqrt{2\sigma}$, after some calculations we obtain the following expression for the integrand in (5.11):
\begin{align*}
x^{2n}e^{-\,\sigma x^{2}} &=\frac{(-1)^{n}}{(2\sigma)^{n}}\sqrt{\frac{2}{\pi}}\int_{0}^{+\infty}H_{2n}(v)e^{-\,v^{2}/2}
\cos{\bigr(xv\sqrt{2\sigma}\bigl)}\d v \\&
=\frac{2}{\sqrt{\pi\sigma}}\,\frac{(-1)^{n}}{(2\sigma)^{n}}\int_{0}^{+\infty}H_{2n}
\biggl(v\sqrt{\frac{2}{\sigma}}\biggr)e^{-\,v^{2}/\sigma}\cos{(2xv)}\d v.
\end{align*}
Therefore, the integral $j(n)$ can be transformed as follows:
\begin{align*}
j(n) &=-\frac{2}{\omega\sqrt{\pi\sigma}}\,\frac{(-1)^{n}}{(2\sigma)^{n}}
\int_{0}^{\,\pi\omega}g\biggl(\frac{x}{\pi\omega}\biggr)
\sin{(2\alpha x)}\int_{0}^{+\infty}H_{2n}\biggl(v\sqrt{\frac{2}{\sigma}}\biggr)e^{-\,v^{2}/\sigma}\cos{(2xv)}\d v\d x\\&
=-\frac{2}{\omega\sqrt{\pi\sigma}}\,\frac{(-1)^{n}}{(2\sigma)^{n}}
\int_{0}^{+\infty}H_{2n}\biggl(v\sqrt{\frac{2}{\sigma}}\biggr)e^{-\,v^{2}/\sigma}
\int_{0}^{\,\pi\omega}g\biggl(\frac{x}{\pi\omega}\biggr)\sin{(2\alpha x)}\cos{(2vx)}\d x\d v\\&
=\frac{-2}{\omega\sqrt{\pi\sigma}}\,\frac{(-1)^{n}}{(2\sigma)^{n}}
\int\limits_{0}^{+\infty}H_{2n}\biggl(v\sqrt{\frac{2}{\sigma}}\biggr)e^{-\,v^{2}/\sigma}
\int\limits_{0}^{\,\pi\omega}g\biggl(\frac{x}{\pi\omega}\biggr)\bigl(\sin{2x(\alpha+v)}\!+\!\sin{2x(\alpha-v)}\bigr)\d x\!\d v\\&
=\frac{\pi}{\omega\sqrt{\pi\sigma}}\,\frac{(-1)^{n}}{(2\sigma)^{n}}
\int\limits_{0}^{+\infty}H_{2n}\biggl(v\sqrt{\frac{2}{\sigma}}\biggr)
e^{-\,v^{2}/\sigma}\int\limits_{0}^{\,\omega}g\biggl(\frac{t}{\omega}\biggr)\bigl(\sin{2\pi t(\alpha+v)}\!+\!\sin{2\pi t(\alpha-v)}\bigr)\d t\d v.
\end{align*}
Using the definition of $F_{\omega}(u)$ and the fact that $F_{\omega}(-u)=-F_{\omega}(u)$, we find that
\begin{multline*}
j(n)\,=\,\frac{\pi}{\sqrt{\pi\sigma}}\,\frac{(-1)^{n}}{(2\sigma)^{n}}
\int_{0}^{+\infty}H_{2n}\biggl(v\sqrt{\frac{2}{\sigma}}\biggr)
e^{-\,v^{2}/\sigma}\bigl(F_{\omega}(v-\alpha)-F_{\omega}(v+\alpha)\bigr)\d v\,=\\
=\,\frac{\pi}{\sqrt{\pi\sigma}}\,\frac{(-1)^{n}}{(2\sigma)^{n}}
\int_{-\infty}^{+\infty}H_{2n}\biggl(v\sqrt{\frac{2}{\sigma}}\biggr)
e^{-\,v^{2}/\sigma}F_{\omega}(v-\alpha)\d v.
\end{multline*}
By Lemma 11,
\begin{multline*}
j(n)\,=\,\frac{\pi}{\sqrt{\pi\sigma}}\,\frac{(-1)^{n}}{(2\sigma)^{n}}
\int_{-\infty}^{+\infty}H_{2n}\biggl(v\sqrt{\frac{2}{\sigma}}\biggr)
e^{-\,v^{2}/\sigma}\bigl(\text{sgn}(v-\alpha)\,+\,\theta K\bigl(\omega(v-\alpha)\bigr)\bigr)\d v\,=\\
=\,\frac{\pi}{\sqrt{2\pi}}\,\frac{(-1)^{n}}{(2\sigma)^{n}}\int_{-\infty}^{+\infty}H_{2n}(u)
e^{-\,u^{2}/2}\bigl(\text{sgn}(u\sqrt{\sigma/2}-\alpha)\,+\,\theta K\bigl(\omega(u\sqrt{\sigma/2}-\alpha)\bigr)\bigr)\d v.
\end{multline*}
Obviously,
\[
\text{sgn}\biggl(u\sqrt{\frac{\sigma}{2}}\,-\,\alpha\biggr)\,=\,
\text{sgn}\biggl(u\,-\,\alpha\sqrt{\frac{2}{\sigma}}\biggr).
\]
Next, by Lemma 9, the contribution to $j(n)$ coming from the $K$-term is estimated as
\[
\ll\,\frac{1}{(2\sigma)^{n}}\,\frac{(n+1)^{1/4}}{\omega\sqrt{\sigma}}\,
\frac{(2n)!}{n!}\biggl(\frac{3}{2}\biggr)^{\!n}\,\ll\,
\frac{(n+1)^{1/4}}{\omega\sqrt{\sigma}}\,\frac{(2n)!}{n!}\biggl(\frac{3}{4\sigma}\biggr)^{\!n}.
\]
In view of (2.5), this term contributes to the sum in (5.10) at most
\begin{multline*}
\ll\frac{H}{\omega\sqrt{\sigma}}\sum\limits_{n = 0}^{\nu}|\Phi_{n}|(n+1)^{1/4}\,\frac{(2n)!}{n!}\biggl(\frac{3}{4\sigma}\biggr)^{\!\!n}\,\ll\,
\frac{H}{\omega\sqrt{\sigma}}\biggl(1\,+\,\sum\limits_{n=2}^{\nu}n^{1/4}
\binom{2n}{n}\biggl(\frac{3}{4\sigma}\biggr)^{\!\!n}
(\log\log{n}+1)^{n}\biggr)\\
\ll\,\frac{H}{\omega\sqrt{\sigma}}\biggl(1\,+\,\sum\limits_{n=2}^{\nu}
\biggl(\frac{3}{\sigma}(\log\log{\nu}+1)\biggr)^{\!\!n}\,\biggr)
\,\ll\,\frac{H}{\omega\sqrt{\sigma}}\biggl(1\,+\,
\sum\limits_{n=2}^{\nu}\biggl(\frac{4}{5}\biggr)^{\!\!n}\,\biggr)\,\ll
\,\frac{H}{\omega\sqrt{\sigma}}.
\end{multline*}
Thus we finally get
\begin{multline*}
I_{2}\,=\,\frac{H}{\pi}\sum\limits_{n=0}^{\nu}(-1)^{n}\Phi_{n}\,
\frac{\pi}{\sqrt{2\pi}}\,\frac{(-1)^{n}}{(2\sigma)^{n}}
\int_{-\infty}^{+\infty}H_{2n}(u)e^{-u^{2}/2}\text{sgn}\bigl(u-\alpha\sqrt{2/\sigma}\bigr)\d u\,+\,O\bigl(H\Delta\bigr)\,=\\
=\,\frac{H}{\sqrt{2\pi}}\sum\limits_{n=0}^{\nu}\frac{\Phi_{n}}{(2\sigma)^{n}}
\int_{-\infty}^{+\infty}H_{2n}(u)e^{-u^{2}/2}\text{sgn}\bigl(u-\alpha\sqrt{2/\sigma}\bigr)\d u\,+\,O\bigl(H\Delta\bigr),
\end{multline*}
where
\[
\Delta\,=\,\frac{1}{\omega\sqrt{\sigma}}\,+\,\frac{e^{-0.5\sigma}}{\sigma}\,+\,\frac{1}{\nu+1}\biggl(\frac{\log\log{(\nu+2)}+1}{\sigma}\biggr)^{\!\nu+1}.
\]
The same formula holds for $I_{1}$ and for the initial integral of the theorem. To end the proof, we have to choose $\omega$.
If
\[
\frac{1}{8\pi\sqrt{e}}\biggl(\frac{\log{H}}{\sigma\log{y}}\biggr)^{\!\!1/2}\,\le\,\frac{1}{4}\,e^{0.5(\log{y})^{1/3}},
\]
then we set
\[
\omega\,=\,\gamma\biggl(\frac{\log{H}}{\sigma\log{y}}\biggr)^{\!\!1/2},\quad \gamma\,=\,\frac{1}{50}\,<\,\frac{1}{8\pi\sqrt{e}}.
\]
In this case,
\[
\frac{1}{\omega\sqrt{\sigma}}\,+\,\frac{e^{-0.5\sigma}}{\sigma}\,\ll\,\biggl(\frac{\log{H}}{\log{y}}\biggr)^{\!\!1/2}
\,+\,\frac{(\log_{2}{y})^{-1}}{\sqrt{\log{y\mathstrut}}}.
\]
Otherwise, if
\[
\frac{1}{4}\,e^{0.5(\log{y})^{1/3}}\,<\,\frac{1}{8\pi\sqrt{e}}\biggl(\frac{\log{H}}{\sigma\log{y}}\biggr)^{\!\!1/2},
\]
then we set
\[
\omega\,=\,\frac{1}{8}\,e^{0.5(\log{y})^{1/3}}.
\]
Thus we have
\begin{align*}
&\frac{1}{\omega\sqrt{\sigma}}\,+\,\frac{e^{-0.5\sigma}}{\sigma}\,\ll\,e^{-0.5(\log{y})^{1/3}}\,+\,
\frac{(\log_{2}{y})^{-1}}{\sqrt{\log{y\mathstrut}}}
\\&
\ll\,\frac{(\log_{2}{y})^{-1}}{\sqrt{\log{y\mathstrut}}}\,\ll\,
\biggl(\frac{\log{H}}{\log{y}}\biggr)^{\!\!1/2}\,+\,\frac{(\log_{2}{y})^{-1}}{\sqrt{\log{y\mathstrut}}}.
\end{align*}
Theorem 3 is proved. $\Box$

Now Corollary 1 follows directly from Theorem 3 and the well-known identity
\[
\chi_{a,b}(u)\,=\,\frac{1}{2}\bigl(\text{sgn}(\beta-u)\,-\,\text{sgn}(u-\alpha)\bigr)\qquad (u\ne\alpha,\beta).
\]

\vspace{0.3cm}

\vspace{0.5cm}

\textbf{6. Proof of Theorem 1}

\vspace{0.5cm}

Suppose that
$0<\varepsilon<10^{-3}$ is an arbitrary small fixed constant, $T\ge T_{0}(\varepsilon)$, $x=T^{\,0.1\varepsilon}$,
$T^{\,c+\varepsilon_{1}}\le H\le T^{\,c+\varepsilon}$, where
$c = \tfrac{27}{82}$, $\varepsilon_{1} =0.9\varepsilon$. Further, let $m=\bigl[\log_{3}{T}\bigr]$, $y = x^{1/(8m+3)}$. Finally, let $V(t) = V_{y}(t)$.
We shall first prove that, for any real $\alpha$,
\[
\int_{T}^{T+H}\text{sgn}\bigl(\pi S(t)-\alpha\bigr)\d t\,=\,-\,\int_{T}^{T+H}\text{sgn}\bigl(V(t)+\alpha\bigr)\d t\,+\,
O\biggl(\frac{H\log_{3}{T}}{\varepsilon\sqrt{\log_{2}{T}}}\biggr),\eqno(6.1)
\]
where the implied constant is absolute. Equation (6.1) shows that the problem of the distribution of the sign of $\pi S(t)-\a$
is transformed into the problem of the distribution of the sign of $V(t)+\a$, and this is handled by Theorem 3.
The idea of proof of (6.1) follows that of Theorem 6.1 from K.-M.~Tsang \cite{Tsang_1984}. Set for brevity $R(t) = \pi S(t)+V(t)$. Given $\alpha$, we have
\[
\pi S(t)-\alpha\,=\,-V(t)-\alpha+R(t).
\]
Denote by $E_{1}$ and $E_{2}$ the sets of $t\in [T,T+H]$ satisfying the conditions
\begin{equation*}
|V(t)+\alpha|\,>\,|R(t)|
\end{equation*}
and
\begin{equation*}
|V(t)+\alpha|\,\le\,|R(t)|,\eqno(6.2)
\end{equation*}
respectively. Since
\[
\text{sgn}\bigl(\pi S(t)-\alpha\bigr)\,=\,\text{sgn}\bigl(-V(t)-\alpha\bigr)\,=\,-\text{sgn}\bigl(V(t)+\alpha\bigr)
\]
for any $t\in E_{1}$, then
\[
\int_{T}^{T+H}\text{sgn}\bigl(\pi S(t)-\alpha\bigr)\d t\,=\,-\int_{T}^{T+H}\text{sgn}\bigl(V(t)+\alpha\bigr)\d t\,+\,2\theta|E_{2}|.
\]
Next, let $\kappa = \pi ec_{0}$, where $c_{0}$ is the constant from Lemma 12, and let $D = \kappa\varepsilon^{-1}\log_{3}{T}$.
Denote by $G_{1}$ and $G_{2}$ the sets of $t\in [T,T+H]$ such that $|V(t)+\alpha|\le D$ and $|R(t)|>D$, respectively.

Suppose that $t\in E_{2}$. If $|R(t)|\le D$ then, by (6.2), we have $|V(t)+\alpha|\le D$ and hence $t\in G_{1}$.
If $|R(t)|>D$ then $t\in G_{2}$. Thus, the set $E_{2}$ is contained in the union of $G_{1}$ and $G_{2}$. Hence,
\[
|E_{2}|\,\le\,|G_{1}|\,+\,|G_{2}|.
\]
On the set $G_{1}$, we have
\[
-D-\alpha\,\le\,V(t)\,\le\,D-\alpha.
\]
Let us define $\xi$ and $ \eta$ by the relations
\[
-D-\alpha\,=\,\xi\sqrt{\frac{\sigma}{2}},\quad D-\alpha\,=\,\eta\sqrt{\frac{\sigma}{2}},\quad\text{where}\quad \sigma\,=\,\sum\limits_{p\le y}\frac{1}{p}.
\]
Setting $\nu = 1$ in (2.8) of Corollary 1, we obtain
\begin{multline*}
|G_{1}|\,\le\,\int_{T}^{T+H}\chi_{(-D-\alpha),D-\alpha}\bigl(V(t)\bigr)\d t\,=\,\frac{H}{\sqrt{2\pi}}\biggl(\;\int_{\xi}^{\eta}\d u\,+\,
O\biggl(\frac{1}{\sigma^{2}}\biggr)\biggr)\,\ll\\
\ll\,H\bigl(\eta-\xi+\sigma^{-2}\bigr)\,\ll\,\frac{HD}{\sqrt{\sigma}}\,\ll\,\frac{H\log_{3}{T}}{\varepsilon\sqrt{\log_{2}{T}}}.
\end{multline*}
To estimate $|G_{2}|$, we use Lemma 12. Indeed,
\[
|G_{2}|D^{2m}\,\le\,\int_{G_{2}}R^{2m}(t)\d t\,\le\,\int_{T}^{T+H}R^{2m}(t)\d t\,\le\,
\frac{c_{0}H}{56\varepsilon}\,(\pi c_{0}\varepsilon^{-1})^{2m},
\]
and hence
\[
|G_{2}|\,\le\,\frac{c_{0}H}{56\varepsilon}\biggl(\frac{\pi c_{0}}{D\varepsilon}\biggr)^{\!2m}\,
\ll\,e^{-2m}\varepsilon^{-1}H\,\ll\,\frac{H}{\varepsilon(\log_{2}{T})^{2}}.
\]
Thus we have
\[
|E_{2}|\,\ll\,\frac{H\log_{3}{T}}{\varepsilon\sqrt{\log_{2}{T\mathstrut}}}.
\]
This proves (6.1)

\vspace{0.5cm}

To prove Theorem 1 it is sufficient to prove the assertion for the case
\begin{equation*}
0.5T^{\,c+\varepsilon}\,\le\,H\,\le\,2T^{\,c+\varepsilon}\qquad(c = \tfrac{27}{82}).\eqno(6.3)
\end{equation*}
Indeed, if $H>2T^{\,c+\varepsilon}$ then we set $h = T^{\,c+\varepsilon} $ and split the interval $T<t\le T+H$ into  segments of the type
\begin{equation*}
T+(s-1)h<t\le T+sh\qquad (s = 1,2,\ldots\,).\eqno(6.4)
\end{equation*}
Each of these intervals, except possibly the last one, has the form $T_{1}<t\le T_{1}+h$, where $T\le T_{1}<2T$. This implies
\[
0.5T_{1}^{\,c+\varepsilon}\,\le\,(0.5T_{1})^{\,c+\varepsilon}\,\le\,T^{\,c+\varepsilon}\,=\,h\,\le\,T_{1}^{\,c+\varepsilon}.
\]
If the length of the last segment from the set (6.4) is less than $h$, we unite it with its left neighbour.
Thus we obtain the segment of the type $T_{1}<t\le T_{1}+h_{1}$, where
\[
0.5T_{1}^{\,c+\varepsilon}\,\le\,h\,<\,h_{1}\,\le\,2h\,\le\,2T_{1}^{\,c+\varepsilon}.
\]
Summation of the both sides of (2.1) over all $T = T_{1}$, $H = h$ (or $H = h_{1}$) leads to the desired assertion.

\medskip
So, let $H$ satisfy (6.3), and let $\alpha<\beta$ be arbitrary real numbers. Using (6.1) we obtain
\begin{multline*}
\frac{1}{2}\int_{T}^{T+H}\bigl(\text{sgn}\bigl(\pi S(t)-\alpha\bigr)\,-\,\text{sgn}\bigl(\pi S(t)-\beta\bigr)\bigr)\d t\,=\\
=\,\frac{1}{2}\int_{T}^{T+H}\bigl(\text{sgn}\bigl(V(t)+\beta\bigr)\,-\,\text{sgn}\bigl(V(t)+\alpha\bigr)\bigr)\d t\,+\,
O\biggl(\frac{H\log_{3}{T}}{\varepsilon\sqrt{\log_{2}{T}}}\biggr),
\end{multline*}
and therefore
\begin{multline*}
\int_{T}^{T+H}\chi_{\alpha,\beta}\bigl(\pi S(t)\bigr)\d t\,=\,\int_{T}^{T+H}\chi_{-\beta,-\alpha}\bigl(V(t)\bigr)\d t\,+\,
O\biggl(\frac{H\log_{3}{T}}{\varepsilon\sqrt{\log_{2}{T}}}\biggr)\,=\\
=\,\frac{H}{\sqrt{2\pi}}\biggl(\;\int_{-\beta\sqrt{2/\sigma}}^{-\alpha\sqrt{2/\sigma}}e^{-v^{2}/2}\d v\,+\,O\bigl(\sigma^{-2}\bigr)\biggr)
\,+\, O\biggl(\frac{H\log_{3}{T}}{\varepsilon\sqrt{\log_{2}{T}}}\biggr)\,=\\
=\,\frac{H}{\sqrt{2\pi}}\;\int_{\alpha\sqrt{2/\sigma}}^{\beta\sqrt{2/\sigma}}e^{-v^{2}/2}\d v\,+
\, O\biggl(\frac{H\log_{3}{T}}{\varepsilon\sqrt{\log_{2}{T}}}\biggr),\qquad(6.5)
\end{multline*}
where
\[
\sigma\,=\,\sum\limits_{p\le y}\frac{1}{p}.
\]
Further, let $a,b$ be given numbers, $a<b$, and let $E = E_{a,b}$ be the set of $t\in [T,T+H]$ satisfying the inequalities
\[
a\,<\,\frac{\pi S(t)\sqrt{2}}{\sqrt{\log\log{T\mathstrut}}}\,\le\,b.
\]
Setting
\[
\alpha\,=\,a\sqrt{\tfrac{1}{2}\log\log{T}},\quad \beta\,=\,b\sqrt{\tfrac{1}{2}\log\log{T}}
\]
in (6.5), we obtain
\begin{equation*}
\text{mes}\bigl\{E_{a,b}\bigr\}\,=\,\frac{H}{\sqrt{2\pi}}\,\biggl(\;\int_{\xi}^{\eta}e^{-v^{2}/2}\d v\,+\,
O\biggl(\frac{\log_{3}{T}}{\varepsilon\sqrt{\log_{2}{T}}}\biggr)\biggr),\eqno(6.6)
\end{equation*}
where
\[
\xi\,=\,\alpha\sqrt{\frac{2}{\sigma}}\,=\,a\sqrt{\frac{1}{\sigma}\,\log\log{T}},\quad
\eta\,=\,\beta\sqrt{\frac{2}{\sigma}}\,=\,b\sqrt{\frac{1}{\sigma}\,\log\log{T}}.
\]
Using the inequalities (3.19) and (3.20) from \cite{Rosser_Schoenfeld_1962}, we get
\begin{multline*}
\sigma\,=\,\log\log{y}\,+\,B\,+\,\frac{\theta}{\log^{2}{y}}\,=\,\log{\biggl(\frac{\varepsilon\log{T}}{10(8m+3)}\biggr)}
\,+\,B\,+\,\frac{\theta}{\log^{2}{y}}\,=\\
=\,\log\log{T}\,-\,\log_{4}{T}\,+\,O\bigl(\log(\varepsilon^{-1})\bigr)
\end{multline*}
and hence
\begin{multline*}
\frac{1}{\sigma}\,\log\log{T}\,=\,\biggl(1\,-\,\frac{\log_{4}{T}+O\bigl(\log(\varepsilon^{-1})\bigr)}{\log_{2}{T}}\biggr)^{\!-1}
\,=\,1\,+\,\frac{\log_{4}{T}}{\log_{2}{T}}\,\bigl(1\,+\,o_{\varepsilon}(1)\bigr),\\
\xi\,=\,a(1+\delta_{1}),\quad \eta\,=\,b(1+\delta_{2}),\quad \delta_{j}\,=\,\frac{\log_{4}{T}}{2\log_{2}{T}}
\,\bigl(1\,+\,o_{\varepsilon}(1)\bigr),\quad j = 1,2.
\end{multline*}
We estimate the error arising after  replacing $\xi, \eta$ with $a$, $b$ in the integral in (6.6).
This error is expressed as $j_{2}-j_{1}$, where
\[
j_{1}\,=\,j_{1}(a)\,=\,\int_{a}^{a(1+\delta_{1})}e^{-v^{2}/2}\d v,\quad j_{2} = j_{2}(b)\,=\,
\int_{b}^{b(1+\delta_{2})}e^{-v^{2}/2}\d v.
\]
Since $|j_{1}(a)| = j_{1}(|a|)$, we have
\[
|j_{1}(a)|\,=\,\int_{|a|}^{|a|(1+\delta_{1})}e^{-v^{2}/2}\d v\,\le\,e^{-a^{2}/2}\int_{|a|}^{|a|(1+\delta_{1})}\d v\,=\,
|a|e^{-a^{2}/2}\delta_{1}\,\le\,\frac{\delta_{1}}{\sqrt{e}}
\]
and, similarly, $|j_{2}(b)|\,\le\,\delta_{2}/\sqrt{e}$. Hence,
\[
|j_{2}-j_{1}|\,\ll\,\delta_{1}+\delta_{2}\,\ll\,\frac{\log_{4}{T}}{\log_{2}{T}}
\]
and
\[
\text{mes}\bigl\{E_{a,b}\bigr\}\,=\,\frac{H}{\sqrt{2\pi}}\,\biggl(\;\int_{a}^{b}e^{-v^{2}/2}\d v\,+\,
O\biggl(\frac{\log_{3}{T}}{\varepsilon\sqrt{\log_{2}{T\mathstrut}}}\biggr)\biggr).
\]
Theorem 1 is proved. $\Box$

\vspace{0.5cm}

\renewcommand{\refname}{\normalsize{References}}

\bigskip
\noindent
Aleksandar P. Ivi\'c\\
Serbian Academy of Sciences and Arts\\
Knez Mihailova 35, 11000 Beograd\\
Serbia\\
{\tt aleksandar.ivic@rgf.bg.ac.rs}

\bigskip
\noindent
Maxim A. Korolev\\
Steklov Mathematical Institute\\
Russian Academy of Sciences\\
119991 Moscow, Gubkina street, 8\\
Russia\\
{\tt korolevma@mi.ras.ru}


\begin{thebibliography}{99}

{\small \bibitem{Bourgain} J.~Bourgain,
{\it Decoupling, exponential sums and the Riemann zeta function,}
J. Amer. Math. Soc., {\bf 30}:1 (2017), 205--224.

\bibitem{Boyarinov_2011} R.N.~Boyarinov, {\it On the value distribution  of the Riemann zeta-function}, Dokl. Akad. Nauk,
{\bf 438}:1(2011), 14--15 (Russian); Doklady Math., \textbf{83}:3 (2011),  290--292 (English).

\bibitem{GonekIvic} S.M.~Gonek and A.~Ivi\'c,
{\it On the distribution of  positive and negative values of Hardy's $Z$-function},
 J.~Number Theory, {\bf 174} (2017), 189--201.

\bibitem{Ivic1} A.~Ivi\'c, {\it The Riemann zeta-function}, John Wiley
\& Sons, New York, 1985, 517 pp.  (2nd. ed.: Dover, Mineola, New York, 2003).

\bibitem{Ivic2} A.~Ivi\'c, {\it The theory of Hardy's $Z$-function}, Cambridge University Press,
Cambridge, 2012, 245 pp.

\bibitem{Karatsuba_1984} A.A.~Karatsuba, {\it On the zeros of the function $\z(s)$ on short intervals of the critical line},
Izv. Akad. Nauk SSSR, Ser. Mat., {\bf 48}:3 (1984), 569--584 (Russian);
Mathematics of the USSR -- Izvestiya, {\bf 24}:3 (1984), 523--537 (English).

\bibitem{Karatsuba_1996} A.A.~Karatsuba, {\it On the function $S(T)$},
Izv. Ross. Akad. Nauk. Ser. Mat., {\bf 60}:5 (1996), 27--56 (Russian);
Izvestiya: Mathematics, {\bf 60}:5 (1996), 901--931 (English).

\bibitem{KarKor_2005} A.A.~Karatsuba and M.A.~Korolev, {\it The argument of the Riemann zeta function},
Uspekhi Mat. Nauk, {\bf 60}:3(363) (2005), 41--96 (Russian); Russian Math. Surveys, {\bf 60}:3 (2005), 433--488 (English).

\bibitem{KarKor_2006} A.A.~Karatsuba and M.A.~Korolev, {\it Behaviour of the argument of the Riemann zeta function on the critical
line}, Uspekhi Mat. Nauk., {\bf 61}:3(369) (2006), 3--92 (Russian); Russian Math. Surveys, {\bf 61}:3 (2006), 389--482 (English).

\bibitem{Korolev_2015}
M.A.~Korolev, {\it Gram's Law in the Theory of Riemann Zeta-Function. Part 1}, Sovrem. Probl. Mat. (Steklov Math. Institute of RAS, Moscow), {\bf 20} (2015), 3--161 (Russian); Proc. Steklov Inst. Math. {\bf 292} (2016), Suppl. iss. 2, 1--146.

\bibitem{Montgomery}  H.L.~Montgomery, {\it The pair correlation of zeros of the
zeta-function},  Analytic number theory (Proc. Symp. Pure Math. Vol. XXIV, St. Louis Univ., St. Louis, Mo., 1972). Amer. Math. Soc., Providence, R.I., 1973, 181--193.

\bibitem{Radziwill_2011} M.~Radziwi{\l}{\l}, {\it Large deviations in Selberg's central limit theorem},
\texttt{arXiv:1108.5092 [math.NT]}, 9 pp.

\bibitem{Rosser_Schoenfeld_1962}
J.B.~Rosser, L.~Schoenfeld, {\it Approximate formulas for some
functions of prime numbers}, Illinois J. Math., \textbf{6}:1 (1962), 64--94.

\bibitem {Selberg_1944} A. Selberg, {\it On the remainder in the formula for $N(T)$, the number of zeros of $\z(s)$ in the
strip $0 < t < T$}, Avh. Norske Videnkaps Akad. Oslo. I. Mat.-Naturv. Klasse, 1944, no. 1, 1--27. (see also: A.~Selberg, Collected Papers. Vol. I. Springer-Verlag, Berlin etc., 1989, 179--204).

\bibitem{Selberg_1946}
A.~Selberg, {\it Contributions to the theory of the Riemann zeta-function}, Arch. Math. Naturvid. \textbf{48}:5 (1946), 89--155
(see also: A.~Selberg, Collected Papers. Vol. I. Springer-Verlag, Berlin etc., 1989, 214--280).

\bibitem{Suetin_2005}
P.K.~Suetin, {\it Classical orthogonal polynomials}. 3rd ed. Fizmatlit, Moscow, 2005, 479 pp. (Russian).

\bibitem{Titchmarsh} E.C.~Titchmarsh, {\it The theory of the Riemann
zeta-function}, 2nd ed., Oxford University Press, Oxford, 1986.

\bibitem{Tsang_1984}
K.-M.~Tsang, {\it The distribution of the values of the Riemann zeta function},
Ph.D. Dissertation, Princeton, 1984 (to be found online at:
{\tt http://www.math.sjsu.  edu/$\sim$goldston/TsangThesis.htm}).

\bibitem{Vaaler_1985}
J.D.~Vaaler, {\it Some extremal functions in Fourier analysis}, Bull. Math. Amer. Soc. (New Ser.), \textbf{12}:2 (1985), 183--216.

}











\end{thebibliography}
\end{document}